\documentclass[12pt]{article}
\def\Dj{\hbox{D\kern-.73em\raise.30ex\hbox{-}
\raise-.30ex\hbox{}}}
\def\dj{\hbox{d\kern-.33em\raise.80ex\hbox{-}
\raise-.80ex\hbox{\kern-.40em}}}
\usepackage{epsfig}
\usepackage{amsmath,amsthm,amsfonts,amssymb,amscd}
\begin{document}

\baselineskip=0.25in

\newtheorem{lem}{Lemma}[section]
\newtheorem{thm}[lem]{Theorem}
\newtheorem{cor}[lem]{Corollary}
\newtheorem*{prop}{Proposition}
\newtheorem{con}[lem]{Conjecture}
\newtheorem{rem}[lem]{Remark}
\newtheorem{defi}[lem]{Definition}
\renewcommand\baselinestretch{1.1}
\def\pf{\noindent {\it Proof.} }
\def\qed{\hfill \rule{4pt}{7pt}}

\begin{center} {\Large \bf All Connected Graphs with Maximum Degree\\
at Most 3 whose Energies are Equal to\\[2mm] the Number of Vertices
\footnotetext[1]{Supported by NSFC No.10831001, PCSIRT and the
``973" program.}
 \footnotetext[2]{Supported by NSFC No.10871166, NSFJS and NSFUJS.}}
 \end{center}

\begin{center}
{ \small  Xueliang Li\footnotemark[1],~~Hongping Ma\footnotemark[2]\\[5pt]
\small Center for Combinatorics and LPMC-TJKLC,\\
\small Nankai University, Tianjin 300071, P.R. China.\\
\small Email: lxl@nankai.edu.cn; mhp@cfc.nankai.edu.cn}
\end{center}

\begin{center}
\small (Received July 8, 2009)
\end{center}

\begin{abstract}
The energy $E(G)$ of a graph $G$ is defined as the sum of the
absolute values of its eigenvalues. Let $S_2$ be the star of order
$2$ (or $K_2$) and $Q$ be the graph obtained from $S_2$ by attaching
two pendent edges to each of the end vertices of $S_2$.
Majstorovi\'{c} et al. conjectured that $S_2$, $Q$ and the complete
bipartite graphs $K_{2,2}$ and $K_{3,3}$ are the only $4$ connected
graphs with maximum degree $\Delta \leq 3$ whose energies are equal
to the number of vertices. This paper is devoted to giving a
confirmative proof to the conjecture.
\end{abstract}

\section{Introduction}

We use Bondy and Murty \cite{BM} for terminology and notations not
defined here. Let $G$ be a simple graph with $n$ vertices and $m$
edges. The {\it cyclomatic number} of a connected graph $G$ is
defined as $c(G)=m-n+1$. A graph $G$ with $c(G)=k$ is called a {\it
$k$-cyclic graph}. In particular, for $c(G)=0,1$ or $2$ we call $G$
a tree, unicyclic or bicyclic graph, respectively. Denote by
$\Delta$ the maximum degree of a graph. The eigenvalues
$\lambda_{1}, \lambda_{2},\ldots, \lambda_{n}$ of the adjacency
matrix $A(G)$ of $G$ are said to be the eigenvalues of the graph
$G$. The $energy$ of $G$ is defined as
$$E=E(G)=\sum_{i=1}^{n}|\lambda_{i}|.$$

For several classes of graphs it has been demonstrated that the
energy exceeds the number of vertices (see, \cite{G1}). In 2007,
Nikiforov \cite{N} showed that for almost all graphs,
$$E=\left(\frac{4}{3\pi}+o(1)\right)n^{3/2}.$$
Thus the number of graphs $G$ satisfying the condition $E(G)<n$ is
relatively small. In \cite{GR}, a connected graph $G$ of order $n$
is called $hypoenergetic$ if $E(G)<n$. For hypoenergetic graphs with
$\Delta \leq 3$, we have the following well known results.

\begin{figure}[ht]
\centering
  \setlength{\unitlength}{0.05 mm}%
  \begin{picture}(1592.3, 330.9)(0,0)
  \put(0,0){\includegraphics{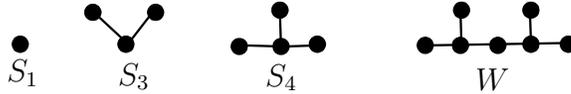}}
  \put(40.00,60.93){\fontsize{11.76}{14.11}\selectfont \makebox(206.6, 82.6)[l]{$S_1$\strut}}
  \put(335.53,52.15){\fontsize{11.76}{14.11}\selectfont \makebox(206.6, 82.6)[l]{$S_3$\strut}}
  \put(727.62,49.22){\fontsize{11.76}{14.11}\selectfont \makebox(206.6, 82.6)[l]{$S_4$\strut}}
  \put(1286.50,37.52){\fontsize{11.76}{14.11}\selectfont \makebox(124.0, 82.6)[l]{$W$\strut}}
  \end{picture}%
 \caption{The hypoenergetic trees with maximum degree at most $3$.} \label{fig1}
\end{figure}

\begin{lem}\cite{GLSZ}\label{lem1.1}
There exist only four hypoenergetic trees with $\Delta \leq 3$,
dipicted in Figure \ref{fig1}.
\end{lem}

\begin{lem}\cite{N2}\label{lem1.2}
Let $G$ be a graph of order $n$ with at least $n$ edges and with no
isolated vertices. If $G$ is quadrangle-free and $\Delta (G)\leq 3$,
then $E(G)>n$.
\end{lem}

The present authors first in \cite{LM1} showed that complete
bipartite graph $K_{2,3}$ is the only hypoenergetic graph among all
unicyclic and bicyclic graphs with $\Delta \leq 3$, and then
recently they obtained the following general result:

\begin{lem}\cite{LM2}\label{lem1.3}
Complete bipartite graph $K_{2,3}$ is the only hypoenergetic
connected cycle-containing (or cyclic) graph with $\Delta \leq 3$.
\end{lem}

Therefore, all connected hypoenergetic graphs with maximum degree at
most 3 have been characterized.

\begin{lem}\cite{LM2}\label{lem1.4}
$S_1, S_3, S_4, W$ and $K_{2,3}$ are the only $5$ hypoenergetic
connected graphs with $\Delta \leq 3$.
\end{lem}

In \cite{MKG} Majstorovi\'{c} et al. proposed the following
conjecture, which is the second half of their Conjecture 3.7.

\begin{figure}[ht]
\centering
  \setlength{\unitlength}{0.05 mm}%
  \begin{picture}(1800.2, 437.5)(0,0)
  \put(0,0){\includegraphics{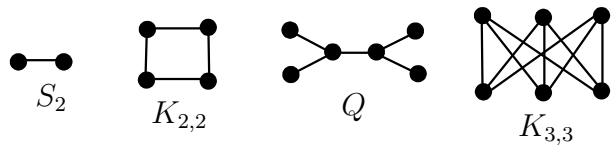}}
  \put(105.89,117.23){\fontsize{11.76}{14.11}\selectfont \makebox(206.6, 82.6)[l]{$S_2$\strut}}
  \put(413.44,71.07){\fontsize{11.76}{14.11}\selectfont \makebox(371.9, 82.6)[l]{$K_{2,2}$\strut}}
  \put(919.87,93.07){\fontsize{11.76}{14.11}\selectfont \makebox(124.0, 82.6)[l]{$Q$\strut}}
  \put(1388.26,37.52){\fontsize{11.76}{14.11}\selectfont \makebox(371.9, 82.6)[l]{$K_{3,3}$\strut}}
  \end{picture}%
  \caption{All connected graphs with maximum degree at most $3$ and $E=n$.} \label{fig2}
\end{figure}

\begin{con}\cite{MKG}\label{con1.5}
There are exactly four connected graphs $G$ with order $n$ and
$\Delta \leq 3$ for which the equality $E(G)=n$ holds, which are
dipicted in Figure \ref{fig2}.
\end{con}

In this paper, we will prove this conjecture.

\section{Main results}

The following results are needed in the sequel.

\begin{lem}\cite{DS}\label{lem2.1}
If $F$ is an edge cut of a graph $G$, then $E(G-F)\leq E(G)$, where
$G-F$ is the subgraph obtained from $G$ by deleting the edges in
$F$.
\end{lem}

\begin{lem}\cite{DS}\label{lem2.2}
Let $F=[S, V\setminus S]$ be an edge cut of a graph $G$ with vertex
set $V$, where $S$ is a nonempty proper subset of $V$. Suppose that
$F$ is not empty and all edges in $F$ are incident to one and only
one vertex in $S$, i.e., the edges in $F$ form a star. Then
$E(G-F)<E(G)$.
\end{lem}

\begin{lem}\cite{BP}\label{lem2.3}
The energy of a graph can not be an odd integer.
\end{lem}

In the following we first show that Conjecture \ref{con1.5} holds
for trees, unicyclic and bicyclic graphs, respectively. Then we show
that Conjecture \ref{con1.5} holds in general.

Let $F$ be an edge cut of a connected graph $F$. If $G-F$ has
exactly two components $G_1$ and $G_2$, then we denote $G-F=G_1+G_2$
for convenience. The following lemma is needed.

\begin{lem}\label{lem2.4}
Let $F$ be an edge cut of a connected graph $G$ of order $n$ such
that $G-F=G_1+G_2$. If $E(G_1)\geq |V(G_1)|$, $E(G_2)\geq |V(G_2)|$
and either at least one of the above inequalities is strict or the
edges in $F$ form a star or both, then $E(G)>n$.
\end{lem}

\pf If $E(G_1)>|V(G_1)|$ or $E(G_2)>|V(G_2)|$, then by Lemma
\ref{lem2.1}, we have
$$E(G)\geq E(G-F)=E(G_1)+E(G_2)>|V(G_1)|+|V(G_2)|=n.$$
Otherwise by Lemma \ref{lem2.2}, we have
$$E(G)>E(G-F)=E(G_1)+E(G_2)\geq |V(G_1)|+|V(G_2)|=n,$$
which completes the proof. \qed

The result Lemma \ref{lem2.4} is easy but useful in our proofs.

\begin{thm}\label{thm2.5}
$S_2$ and $Q$ are the only two trees $T$ with order $n$ and $\Delta
\leq 3$ for which the equality $E(T)=n$ holds.
\end{thm}

\pf Let $T$ be a tree with $n$ vertices and $\Delta \leq 3$. From
Table 2 of \cite{CDS}, we know that $S_2$ and $Q$ are the only two
trees with $\Delta \leq 3$ and $n\leq 10$ for which the equality
$E=n$ holds. By Lemma \ref{lem2.3}, we may assume that $n\geq 12$ is
even. We will prove that $E(T)>n$.

We divide the trees with $\Delta \leq 3$ into two classes: {\bf
Class 1} contains the trees $T$ that have an edge $e$, such that
$T-e=T_1+T_2$ and $T_1, T_2\not\cong S_1, S_3, S_4, W$. {\bf Class
2} contains the trees $T$ in which there exists no edge $e$, such
that $T-e=T_1+T_2$ and $T_1, T_2\not\cong S_1, S_3, S_4, W$, i.e.,
for any edge $e$ of $T$ at least one of components of $T-e$ is
isomorphic to a tree in $\{S_1, S_3, S_4, W\}$.

\noindent {\bf Case 1.} $T$ belongs to Class 1. Then there exists an
edge $e$ such that $T-e=T_1+T_2$ and $T_1, T_2\not\cong S_1, S_3,
S_4, W$. Hence by Lemmas \ref{lem1.1} and \ref{lem2.2}, we have
$E(T)>E(T-e)=E(T_1)+E(T_2)\geq |V(T_1)|+|V(T_2)|=n$, which completes
the proof.

\noindent {\bf Case 2.} $T$ belongs to Class 2. Consider the center
of $T$. There are two subcases: either $T$ has a (unique) center
edge $e$ or a (unique) center vertex $v$.

\noindent {\bf Subcase 2.1.} $T$ has a center edge $e$. The two
fragments attached to $e$ will be denoted by $T_1$ and $T_2$, i.e.,
$T-e=T_1+T_2$.

Without loss of generality, we assume that $T_1$ is isomorphic to a
tree in $\{S_1, S_3, S_4$, $W\}$.

If $T_1$ is isomorphic to a tree in $\{S_1, S_3, S_4\}$, then it is
easy to see that $n\leq 11$, which is a contradiction.

\begin{figure}[ht]
\centering
  \setlength{\unitlength}{0.05 mm}%
  \begin{picture}(1993.4, 2404.1)(0,0)
  \put(0,0){\includegraphics{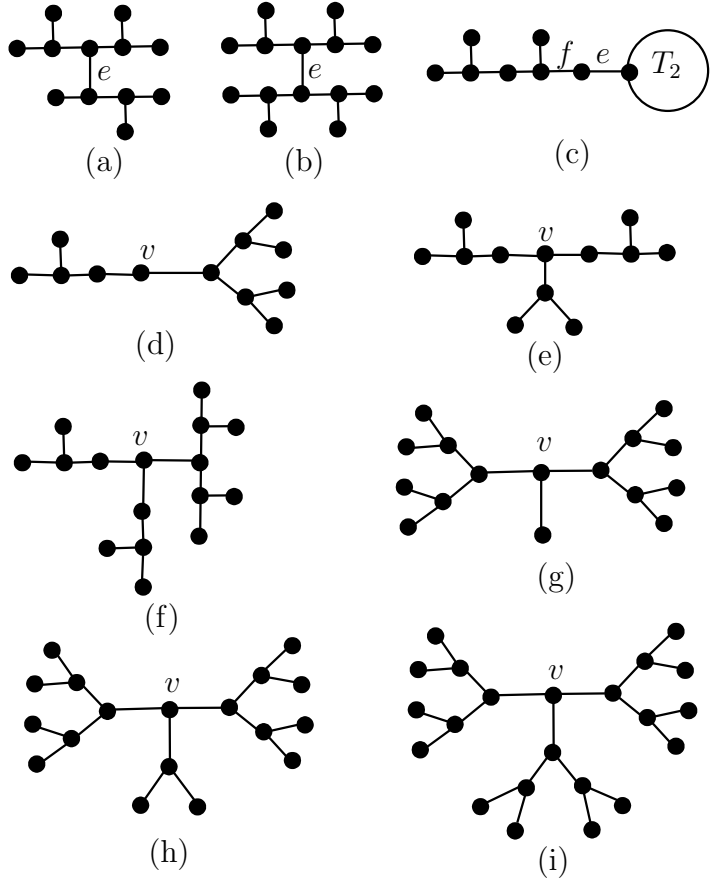}}
  \put(833.42,2146.42){\fontsize{11.76}{14.11}\selectfont \makebox(124.0, 82.6)[l]{$e$\strut}}
  \put(273.06,2139.18){\fontsize{11.76}{14.11}\selectfont \makebox(124.0, 82.6)[l]{$e$\strut}}
  \put(239.24,1897.64){\fontsize{11.76}{14.11}\selectfont \makebox(124.0, 82.6)[l]{(a)\strut}}
  \put(768.21,1897.64){\fontsize{11.76}{14.11}\selectfont \makebox(124.0, 82.6)[l]{(b)\strut}}
  \put(1599.09,2185.07){\fontsize{11.76}{14.11}\selectfont \makebox(124.0, 82.6)[l]{$e$\strut}}
  \put(1746.76,2158.63){\fontsize{11.76}{14.11}\selectfont \makebox(206.6, 82.6)[l]{$T_2$\strut}}
  \put(1486.12,1920.44){\fontsize{11.76}{14.11}\selectfont \makebox(124.0, 82.6)[l]{(c)\strut}}
  \put(1492.61,2186.23){\fontsize{11.76}{14.11}\selectfont \makebox(124.0, 82.6)[l]{$f$\strut}}
  \put(1445.24,1713.09){\fontsize{11.76}{14.11}\selectfont \makebox(124.0, 82.6)[l]{$v$\strut}}
  \put(386.34,1663.54){\fontsize{11.76}{14.11}\selectfont \makebox(124.0, 82.6)[l]{$v$\strut}}
  \put(366.56,1169.11){\fontsize{11.76}{14.11}\selectfont \makebox(124.0, 82.6)[l]{$v$\strut}}
  \put(1441.29,1152.01){\fontsize{11.76}{14.11}\selectfont \makebox(124.0, 82.6)[l]{$v$\strut}}
  \put(452.42,511.23){\fontsize{11.76}{14.11}\selectfont \makebox(124.0, 82.6)[l]{$v$\strut}}
  \put(1472.93,546.83){\fontsize{11.76}{14.11}\selectfont \makebox(124.0, 82.6)[l]{$v$\strut}}
  \put(373.73,1411.12){\fontsize{11.76}{14.11}\selectfont \makebox(124.0, 82.6)[l]{(d)\strut}}
  \put(1413.60,1381.63){\fontsize{11.76}{14.11}\selectfont \makebox(124.0, 82.6)[l]{(e)\strut}}
  \put(397.46,683.32){\fontsize{11.76}{14.11}\selectfont \makebox(124.0, 82.6)[l]{(f)\strut}}
  \put(1441.29,792.66){\fontsize{11.76}{14.11}\selectfont \makebox(124.0, 82.6)[l]{(g)\strut}}
  \put(409.33,58.36){\fontsize{11.76}{14.11}\selectfont \makebox(124.0, 82.6)[l]{(h)\strut}}
  \put(1441.72,37.52){\fontsize{11.76}{14.11}\selectfont \makebox(124.0, 82.6)[l]{(i)\strut}}
  \end{picture}%
  \caption{The graphs in the proof of Theorem \ref{thm2.5}. }\label{fig3}
\end{figure}

If $T_1\cong W$ and it is attached to the center edge $e$ through
the vertex of degree 2, then it is easy to see that $T$ must be the
tree as given in Figure \ref{fig3} (a) or (b). By direct computing,
we have that $E(T)=12.61708>12=n$ in the former case while
$E(T)=14.91128>14=n$ in the latter case. If $T_1\cong W$ and it is
attached to the center edge $e$ through a pendent vertex, see Figure
\ref{fig3} (c). Since $T$ belongs to Class 2, deleting the edge $f$,
we then have that $T_2\cup e$ is isomorphic to a tree in $\{S_1,
S_3, S_4, W\}$, which contradicts to the fact that $e$ is the center
edge of $T$.

\noindent {\bf Subcase 2.2.} $T$ has a center vertex $v$. If $v$ is
of degree 2, then the two fragments attached to it will be denoted
by $T_1$ and $T_2$. If $v$ is of degree 3, then the three fragments
attached to it will be denoted by $T_1$, $T_2$ and $T_3$.

Let $v_i$ be the adjacent vertex of $v$ in $T_i$. Denote
$T-vv_1=T_1+T'_2$. Since $T$ belongs to Class 2, either $T_1$ or
$T'_2$ is isomorphic to a tree in $\{S_1, S_3, S_4, W\}$.

\noindent {\bf Subsubcase 2.2.1.} $T'_2$ is isomorphic to a tree in
$\{S_1, S_3, S_4, W\}$.

Clearly $T'_2\not \cong S_1$. If $T'_2\cong S_3$ or $S_4$, then it
is easy to see that $n\leq 7$, which is a contradiction. If
$T'_2\cong W$ and $v$ is of degree $3$, then it is easy to see that
$n\leq 10$, which is a contradiction. If $T'_2\cong W$ and $v$ is of
degree $2$, i.e., $N(v)=\{v_1,v_2\}$. Consider $T-vv_2$, since $T$
belongs to Class 2, we have that $T_1\cup vv_1$ is isomorphic to a
tree in $\{S_1, S_3, S_4, W\}$. By the fact that $v$ is the center
of $T$, we have that $T_1\cup vv_1\cong W$, and so $n=13$, which is
a contradiction.

\noindent {\bf Subsubcase 2.2.2.} $T_1$ is isomorphic to a tree in
$\{S_1, S_3, S_4, W\}$.

If $T_1\cong S_1$, then it is easy to see that $n\leq 4$, which is a
contradiction.

If $T_1\cong S_3$ and $v_1$ is of degree 2 in $T_1$, then it is easy
to see that $n\leq 10$, which is a contradiction. If $T_1\cong S_3$
and $v_1$ is a pendent vertex in $T_1$, denote by $u$ the unique
adjacent vertex of $v_1$ in $T_1$. Since $T$ belongs to Class 2,
deleting the edge $uv_1$, we then have that $T'_2\cup vv_1$ is
isomorphic to a tree in $\{S_1, S_3, S_4, W\}$, and so $n\leq 9$,
which is a contradiction.

If $T_1\cong S_4$ or $T_1\cong W$ and $v_1$ is of degree 2 in $T_1$,
then by the facts that $T$ belongs to Class 2, $v$ is the center of
$T$ and $n$ is even, it is not hard to obtain that $T_2$, $T_3$ must
be isomorphic to a tree in $\{S_1, S_3, S_4, W\}$, and at least one
of $T_2$ and $T_3$ is isomorphic to a tree in $\{S_4, W\}$, and if
$T_2$ ($T_3$, respectively) is isomorphic to $W$, then $v_2$ ($v_3$,
respectively) is of degree 2 in $T_2$ ($T_3$, respectively). Hence
there are 6 such trees, as given in Figure \ref{fig3} (d), (e), (f),
(g), (h) and (i). The energy of these trees are 12.72729 ($>12=n$),
12.65406 ($>12=n$), 16.81987 ($>16=n$), 16.77215 ($>16=n$), 19.18674
($>18=n$) and 23.38426 ($>22=n$), respectively.

If $T_1\cong W$ and $v_1$ is a pendent vertex in $T_1$, denote by
$u$ the unique adjacent vertex of $v_1$ in $T_1$. Since $T$ belongs
to Class 2, deleting the edge $uv_1$, we then have that $T'_2\cup
vv_1$ is isomorphic to a tree in $\{S_1, S_3, S_4, W\}$, which
contradicts to the fact that $v$ is the center vertex of $T$. The
proof is thus complete. \qed

From Table 1 of \cite{CDS}, we know that $K_{2,2}$ is the only
connected graph of order 4 with $\Delta \leq 3$ and $E=4$. From
Tables 1 and 2 of \cite{CP}, we know that $K_{3,3}$ is the only
connected cycle-containing graph of order 6 with $\Delta \leq 3$ and
$E=6$.

\begin{thm}\label{thm2.6}
$K_{2,2}$ is the only unicyclic graph with $\Delta \leq 3$ for which
the equality $E=n$ holds.
\end{thm}

\pf Let $G\not\cong K_{2,2}$ be a unicyclic graph of order $n$ with
$\Delta \leq 3$. It is sufficient to show that $E(G)>n$. By Lemmas
\ref{lem1.2} and \ref{lem2.3}, we can assume that $n\geq 8$ is even
and $G$ contains a quadrangle $C=x_1x_2x_3x_4x_1$. We distinguish
the following four cases:

\noindent {\bf Case 1.} There exists an edge $e$ on $C$ such that
the end vertices of $e$ are of degree 2.

Without loss of generality, we assume that $d(x_1)=d(x_4)=2$. Let
$F=\{x_1x_2,x_4x_3\}$, then $G-F=G_1+G_2$, where $G_1\cong S_2$ and
$G_2$ is a tree of order at least $6$ since $n\geq 8$. Since
$\Delta(G)\leq 3$, $G_2$ can not be isomorphic to $W$ or $Q$.
Therefore we have $E(G_1)=|V(G_1)|$ and $E(G_2)>|V(G_2)|$ by Lemma
\ref{lem1.1} and Theorem \ref{thm2.5}. It follows from Lemma
\ref{lem2.4} that $E(G)>n$.

\noindent {\bf Case 2.} There exist exactly two nonadjacent vertices
$x_i$ and $x_j$ on $C$ such that $d(x_i)=d(x_j)=2$.

Without loss of generality, we assume that $d(x_2)=d(x_4)=2$,
$d(x_1)=d(x_3)=3$. Let $y_3$ be the adjacent vertex of $x_3$ outside
$C$. Then $G-x_3y_3=G_1+G_2$, where $G_1$ is a unicyclic graph and
$G_2$ is a tree. Notice that $E(G_1)\geq |V(G_1)|$ by Lemma
\ref{lem1.3}. If $G_2\not\cong S_1,S_3,S_4,W$, then we have
$E(G_2)\geq |V(G_2)|$ by Lemma \ref{lem1.1} and so $E(G)>
E(G-x_3y_3)\geq n$ by Lemma \ref{lem2.4}. Therefore we only need to
consider the following four subcases.

\noindent {\bf Subcase 2.1.} $G_2\cong S_1$. Let
$F=\{x_2x_3,x_3x_4\}$, then $G-F=G'_1+G'_2$, where $G'_2\cong S_2$
and $G'_1$ is a tree of order at least $6$ since $n\geq 8$. If
$G'_1\cong W$, then $n=9$, which is a contradiction. Otherwise, it
follows from Lemmas \ref{lem1.1} and \ref{lem2.4} that $E(G)>n$.

\begin{figure}[ht]
\centering
 \setlength{\unitlength}{0.05 mm}%
  \begin{picture}(2674.2, 2001.3)(0,0)
  \put(0,0){\includegraphics{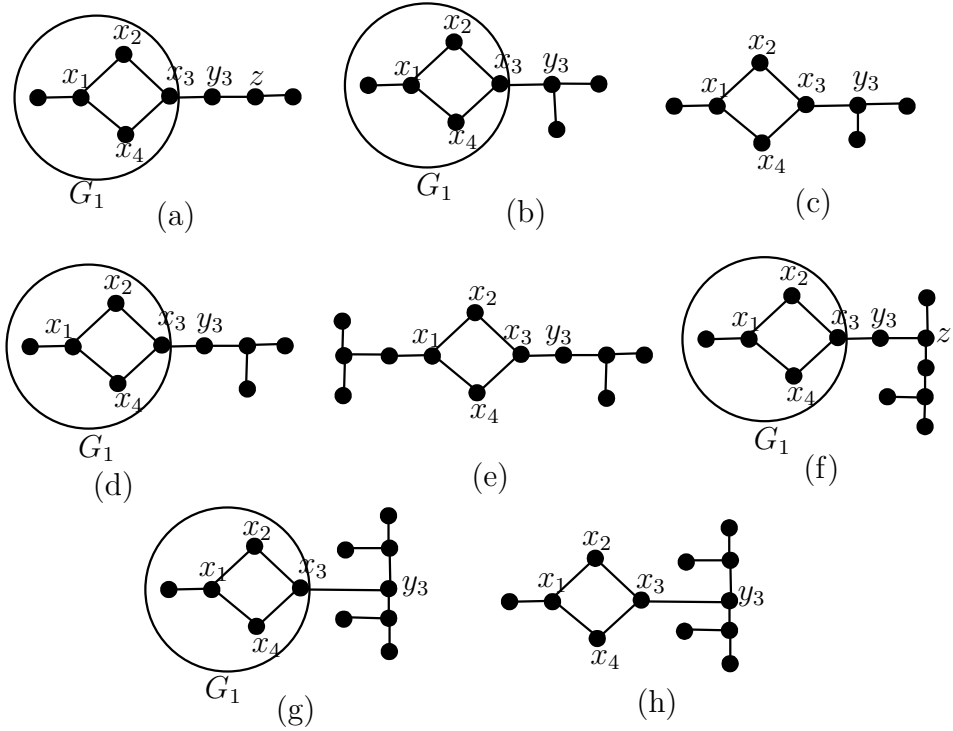}}
  \put(436.71,1351.11){\fontsize{11.76}{14.11}\selectfont \makebox(124.0, 82.6)[l]{(a)\strut}}
  \put(322.97,1845.12){\fontsize{11.76}{14.11}\selectfont \makebox(206.6, 82.6)[l]{$x_2$\strut}}
  \put(328.61,1531.95){\fontsize{11.76}{14.11}\selectfont \makebox(206.6, 82.6)[l]{$x_4$\strut}}
  \put(187.54,1729.45){\fontsize{11.76}{14.11}\selectfont \makebox(206.6, 82.6)[l]{$x_1$\strut}}
  \put(464.04,1737.91){\fontsize{11.76}{14.11}\selectfont \makebox(206.6, 82.6)[l]{$x_3$\strut}}
  \put(571.25,1737.91){\fontsize{11.76}{14.11}\selectfont \makebox(206.6, 82.6)[l]{$y_3$\strut}}
  \put(683.94,1723.57){\fontsize{11.76}{14.11}\selectfont \makebox(124.0, 82.6)[l]{$z$\strut}}
  \put(1365.89,1358.30){\fontsize{11.76}{14.11}\selectfont \makebox(124.0, 82.6)[l]{(b)\strut}}
  \put(206.02,1410.19){\fontsize{11.76}{14.11}\selectfont \makebox(206.6, 82.6)[l]{$G_1$\strut}}
  \put(1200.52,1874.46){\fontsize{11.76}{14.11}\selectfont \makebox(206.6, 82.6)[l]{$x_2$\strut}}
  \put(1214.24,1566.29){\fontsize{11.76}{14.11}\selectfont \makebox(206.6, 82.6)[l]{$x_4$\strut}}
  \put(1069.44,1749.90){\fontsize{11.76}{14.11}\selectfont \makebox(206.6, 82.6)[l]{$x_1$\strut}}
  \put(1335.10,1769.01){\fontsize{11.76}{14.11}\selectfont \makebox(206.6, 82.6)[l]{$x_3$\strut}}
  \put(1464.90,1771.83){\fontsize{11.76}{14.11}\selectfont \makebox(206.6, 82.6)[l]{$y_3$\strut}}
  \put(2134.36,1385.03){\fontsize{11.76}{14.11}\selectfont \makebox(124.0, 82.6)[l]{(c)\strut}}
  \put(1129.52,1430.57){\fontsize{11.76}{14.11}\selectfont \makebox(206.6, 82.6)[l]{$G_1$\strut}}
  \put(269.37,631.98){\fontsize{11.76}{14.11}\selectfont \makebox(124.0, 82.6)[l]{(d)\strut}}
  \put(2010.96,1819.07){\fontsize{11.76}{14.11}\selectfont \makebox(206.6, 82.6)[l]{$x_2$\strut}}
  \put(2035.48,1501.81){\fontsize{11.76}{14.11}\selectfont \makebox(206.6, 82.6)[l]{$x_4$\strut}}
  \put(1885.94,1704.95){\fontsize{11.76}{14.11}\selectfont \makebox(206.6, 82.6)[l]{$x_1$\strut}}
  \put(2142.12,1717.79){\fontsize{11.76}{14.11}\selectfont \makebox(206.6, 82.6)[l]{$x_3$\strut}}
  \put(2289.30,1724.42){\fontsize{11.76}{14.11}\selectfont \makebox(206.6, 82.6)[l]{$y_3$\strut}}
  \put(2095.76,1204.77){\fontsize{11.76}{14.11}\selectfont \makebox(206.6, 82.6)[l]{$x_2$\strut}}
  \put(2114.12,884.74){\fontsize{11.76}{14.11}\selectfont \makebox(206.6, 82.6)[l]{$x_4$\strut}}
  \put(1967.52,1083.03){\fontsize{11.76}{14.11}\selectfont \makebox(206.6, 82.6)[l]{$x_1$\strut}}
  \put(2231.86,1088.19){\fontsize{11.76}{14.11}\selectfont \makebox(206.6, 82.6)[l]{$x_3$\strut}}
  \put(2336.66,1092.53){\fontsize{11.76}{14.11}\selectfont \makebox(206.6, 82.6)[l]{$y_3$\strut}}
  \put(2510.26,1040.18){\fontsize{11.76}{14.11}\selectfont \makebox(124.0, 82.6)[l]{$z$\strut}}
  \put(1714.73,54.45){\fontsize{11.76}{14.11}\selectfont \makebox(124.0, 82.6)[l]{(h)\strut}}
  \put(2027.58,755.24){\fontsize{11.76}{14.11}\selectfont \makebox(206.6, 82.6)[l]{$G_1$\strut}}
  \put(667.18,527.12){\fontsize{11.76}{14.11}\selectfont \makebox(206.6, 82.6)[l]{$x_2$\strut}}
  \put(692.26,214.69){\fontsize{11.76}{14.11}\selectfont \makebox(206.6, 82.6)[l]{$x_4$\strut}}
  \put(551.20,417.83){\fontsize{11.76}{14.11}\selectfont \makebox(206.6, 82.6)[l]{$x_1$\strut}}
  \put(815.84,425.02){\fontsize{11.76}{14.11}\selectfont \makebox(206.6, 82.6)[l]{$x_3$\strut}}
  \put(1091.42,380.33){\fontsize{11.76}{14.11}\selectfont \makebox(206.6, 82.6)[l]{$y_3$\strut}}
  \put(567.06,93.31){\fontsize{11.76}{14.11}\selectfont \makebox(206.6, 82.6)[l]{$G_1$\strut}}
  \put(1570.40,499.53){\fontsize{11.76}{14.11}\selectfont \makebox(206.6, 82.6)[l]{$x_2$\strut}}
  \put(1592.96,183.81){\fontsize{11.76}{14.11}\selectfont \makebox(206.6, 82.6)[l]{$x_4$\strut}}
  \put(1452.00,388.51){\fontsize{11.76}{14.11}\selectfont \makebox(206.6, 82.6)[l]{$x_1$\strut}}
  \put(1714.28,384.98){\fontsize{11.76}{14.11}\selectfont \makebox(206.6, 82.6)[l]{$x_3$\strut}}
  \put(1984.84,353.22){\fontsize{11.76}{14.11}\selectfont \makebox(206.6, 82.6)[l]{$y_3$\strut}}
  \put(298.30,1180.74){\fontsize{11.76}{14.11}\selectfont \makebox(206.6, 82.6)[l]{$x_2$\strut}}
  \put(326.52,867.85){\fontsize{11.76}{14.11}\selectfont \makebox(206.6, 82.6)[l]{$x_4$\strut}}
  \put(141.76,1064.07){\fontsize{11.76}{14.11}\selectfont \makebox(206.6, 82.6)[l]{$x_1$\strut}}
  \put(441.47,1085.16){\fontsize{11.76}{14.11}\selectfont \makebox(206.6, 82.6)[l]{$x_3$\strut}}
  \put(545.21,1081.17){\fontsize{11.76}{14.11}\selectfont \makebox(206.6, 82.6)[l]{$y_3$\strut}}
  \put(1275.83,655.60){\fontsize{11.76}{14.11}\selectfont \makebox(124.0, 82.6)[l]{(e)\strut}}
  \put(230.28,735.10){\fontsize{11.76}{14.11}\selectfont \makebox(206.6, 82.6)[l]{$G_1$\strut}}
  \put(2159.44,679.23){\fontsize{11.76}{14.11}\selectfont \makebox(124.0, 82.6)[l]{(f)\strut}}
  \put(1268.09,1151.64){\fontsize{11.76}{14.11}\selectfont \makebox(206.6, 82.6)[l]{$x_2$\strut}}
  \put(1271.80,838.08){\fontsize{11.76}{14.11}\selectfont \makebox(206.6, 82.6)[l]{$x_4$\strut}}
  \put(1117.70,1038.86){\fontsize{11.76}{14.11}\selectfont \makebox(206.6, 82.6)[l]{$x_1$\strut}}
  \put(1361.10,1046.84){\fontsize{11.76}{14.11}\selectfont \makebox(206.6, 82.6)[l]{$x_3$\strut}}
  \put(1472.82,1046.84){\fontsize{11.76}{14.11}\selectfont \makebox(206.6, 82.6)[l]{$y_3$\strut}}
  \put(756.64,37.52){\fontsize{11.76}{14.11}\selectfont \makebox(124.0, 82.6)[l]{(g)\strut}}
  \end{picture}%
  \caption{The graphs in the proof of Theorem \ref{thm2.6}. }\label{fig4}
\end{figure}

\noindent {\bf Subcase 2.2.} $G_2\cong S_3$. Then $G$ must have the
structure as given in Figure \ref{fig4} (a) or (b). In the former
case, $G-y_3z=G'_1+G'_2$, where $G'_1$ is a unicyclic graph and
$G'_2\cong S_2$. It follows from Lemmas \ref{lem1.4} and
\ref{lem2.4} that $E(G)>n$. In the latter case,
$G-\{x_1x_2,x_4x_3\}=G'_1+G'_2$, where $G'_2$ is the tree of order
$5$ containing $x_3$ and $G'_1$ is a tree of order at least $3$. By
Lemma \ref{lem1.1} and Theorem \ref{thm2.5}, we have
$E(G'_2)>|V(G'_2)|$. If $G'_1\not\cong S_3,S_4,W$, then we have
$E(G)>n$ by Lemmas \ref{lem1.1} and \ref{lem2.4}. Since
$\Delta(G)\leq 3$, $G'_1$ can not be isomorphic to $S_4$ or $W$. If
$G'_1\cong S_3$, then $G$ must be the graph as given in Figure
\ref{fig4} (c). By choosing the edge cut $\{x_1x_2,x_1x_4\}$, we can
similarly obtain that $E(G)>n$.

\noindent {\bf Subcase 2.3.} $G_2\cong S_4$. Then $G$ must have the
structure as given in Figure \ref{fig4} (d). Let
$F=\{x_2x_3,x_3x_4\}$, then $G-F=G'_1+G'_2$, where $G'_2$ is the
tree of order $5$ containing $x_3$ and $G'_1$ is a tree of order at
least $4$. By Lemma \ref{lem1.1} and Theorem \ref{thm2.5}, we have
$E(G'_2)>|V(G'_2)|$. If $G'_1\not\cong S_4,W$, then we have $E(G)>n$
by Lemmas \ref{lem1.1} and \ref{lem2.4}. If $G'_1\cong S_4$, then
$n=9$, which is a contradiction. If $G'_1\cong W$, then $G$ must be
the graph as given in Figure \ref{fig4} (e). By choosing the edge
cut $\{x_1x_2,x_3x_4\}$, we can similarly obtain that $E(G)>n$.

\noindent {\bf Subcase 2.4.} $G_2\cong W$. Then $G$ must have the
structure as given in Figure \ref{fig4} (f) or (g). In the former
case, $G-y_3z=G'_1+G'_2$, where $G'_1$ is a unicyclic graph and
$G'_2$ is a tree of order $6$. It follows from Lemmas \ref{lem1.4}
and \ref{lem2.4} that $E(G)>n$. In the latter case,
$G-\{x_2x_3,x_3x_4\}=G'_1+G'_2$, where $G'_2$ is the tree of order
$8$ containing $x_3$ and $G'_1$ is a tree of order at least $4$. If
$G'_1\not\cong S_4,W$, then we have $E(G)>n$ by Lemmas \ref{lem1.1}
and \ref{lem2.4}. If $G'_1\cong S_4$, then $G$ must be the graph as
given in Figure \ref{fig4} (h). By choosing the edge cut
$\{x_1x_2,x_1x_4\}$, we can similarly obtain that $E(G)>n$. If
$G'_1\cong W$, then $n=15$, which is a contradiction.

\noindent {\bf Case 3.} There exists exactly one vertices $x_i$ on
$C$ such that $d(x_i)=2$.

Without loss of generality, we assume that $d(x_1)=2$. Let
$F=\{x_1x_4,x_2x_3\}$, then $G-F=G_1+G_2$, where $G_1$ is the tree
of order at least $3$ containing  $x_1$ and $G_2$ is a tree of order
at least $4$. Since $\Delta(G)\leq 3$, $G_1, G_2$ can not be
isomorphic to $S_4$, $W$ or $Q$. If $G_1\not\cong S_3$, then we have
$E(G)>n$ by Lemmas \ref{lem1.1}, \ref{lem2.4} and Theorem
\ref{thm2.5}. If $G_1\cong S_3$, then
$G-\{x_1x_2,x_2x_3\}=G'_1+G'_2$, where $G'_1$ is the tree of order
at least $5$ containing $x_1$ and $G'_2\cong S_2$. If $G'_1\not\cong
W$, then we have $E(G)>n$ by Lemmas \ref{lem1.1} and \ref{lem2.4}.
If $G'_1\cong W$, then $n=9$, which is a contradiction.

\noindent {\bf Case 4.} $d(x_1)=d(x_2)=d(x_3)=d(x_4)=3$.

Let $F=\{x_1x_4,x_2x_3\}$, then $G-F=G_1+G_2$, where $G_1$ and $G_2$
are trees of order at least $4$ and it is easy to see that $G_1,
G_2$ can not be isomorphic to $S_4$, $W$ or $Q$. So it follows from
Lemmas \ref{lem1.1}, \ref{lem2.4} and Theorem \ref{thm2.5} that
$E(G)>n$. The proof is thus complete. \qed

\begin{thm}\label{thm2.7}
There does not exist any bicyclic graph with $\Delta \leq 3$ for
which the equality $E=n$ holds.
\end{thm}

\pf Let $G$ be a bicyclic graph of order $n$ with $\Delta \leq 3$.
We know that $E(G)\neq n$ for $n=4$ or $6$. By Lemmas \ref{lem1.2}
and \ref{lem2.3}, we may assume that $n\geq 8$ is even and $G$
contains a quadrangle. Then we will show that $E(G)> n$.

If the cycles in $G$ are disjoint, then it is clear that there
exists a path $P$ connecting the two cycles in $G$. For any edge $e$
on $P$, we have $G-e=G_1+G_2$, where $G_1$ and $G_2$ are unicyclic
graphs. By Lemma \ref{lem1.3}, we have $E(G_1)\geq |V(G_1)|$ and
$E(G_2)\geq |V(G_2)|$. Therefore we have $E(G)>n$ by Lemma
\ref{lem2.4}. Otherwise, the cycles in $G$ have two or more common
vertices. Then we can assume that $G$ contains a subgraph as given
in Figure \ref{fig5} (a), where $P_1,P_2,P_3$ are paths in $G$. We
distinguish the following three cases:
\begin{figure}
\centering
 \setlength{\unitlength}{0.05 mm}%
  \begin{picture}(2805.0, 4124.6)(0,0)
  \put(0,0){\includegraphics{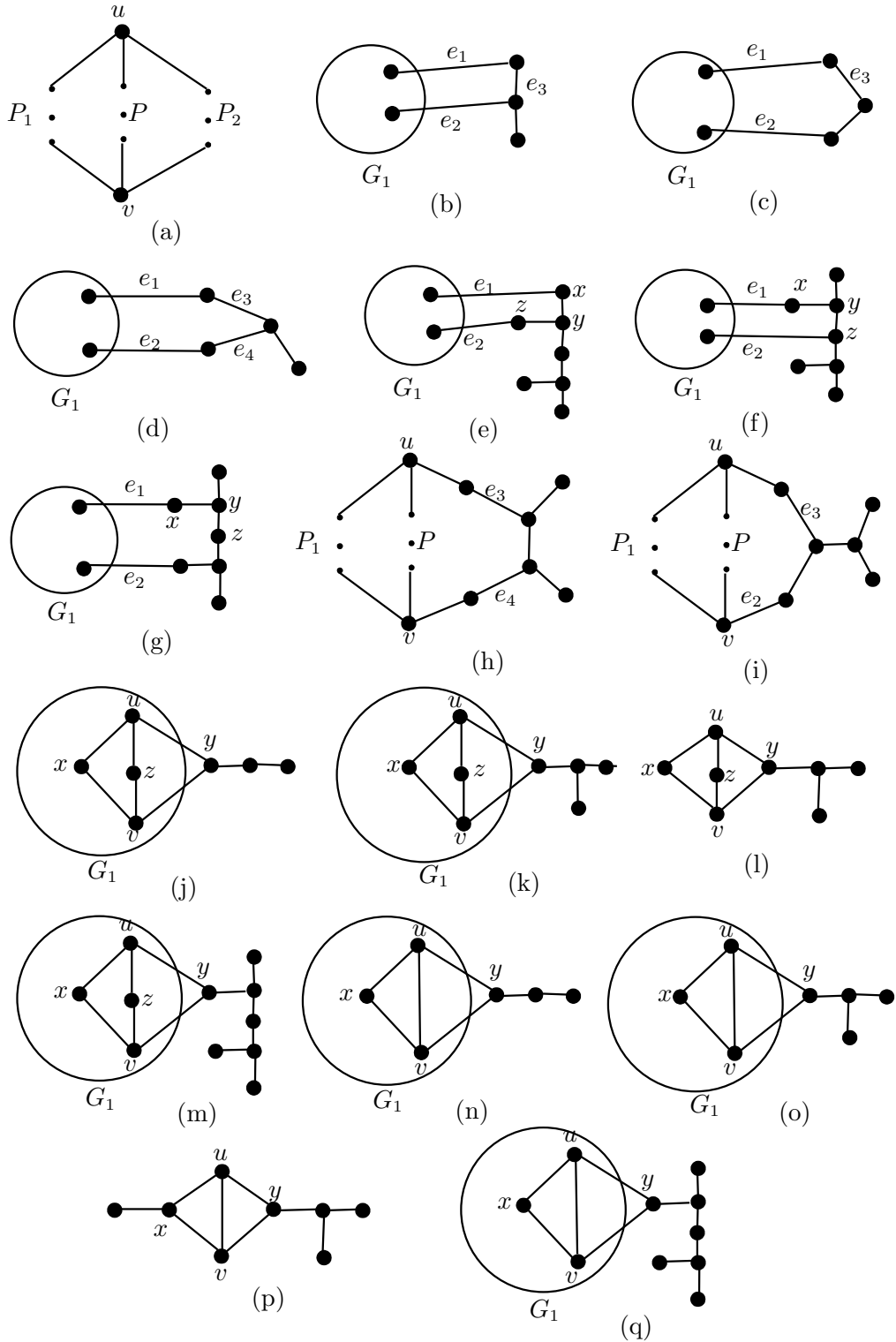}}
  \put(404.40,3692.94){\fontsize{10.69}{12.82}\selectfont \makebox(112.7, 75.1)[l]{$P$\strut}}
  \put(667.92,3693.14){\fontsize{10.69}{12.82}\selectfont \makebox(187.8, 75.1)[l]{$P_2$\strut}}
  \put(40.00,3693.14){\fontsize{10.69}{12.82}\selectfont \makebox(187.8, 75.1)[l]{$P_1$\strut}}
  \put(471.26,3337.36){\fontsize{10.69}{12.82}\selectfont \makebox(112.7, 75.1)[l]{(a)\strut}}
  \put(1312.48,3418.32){\fontsize{10.69}{12.82}\selectfont \makebox(112.7, 75.1)[l]{(b)\strut}}
  \put(2272.66,3427.17){\fontsize{10.69}{12.82}\selectfont \makebox(112.7, 75.1)[l]{(c)\strut}}
  \put(428.44,2742.98){\fontsize{10.69}{12.82}\selectfont \makebox(112.7, 75.1)[l]{(d)\strut}}
  \put(1428.46,2738.99){\fontsize{10.69}{12.82}\selectfont \makebox(112.7, 75.1)[l]{(e)\strut}}
  \put(1445.81,3182.95){\fontsize{10.69}{12.82}\selectfont \makebox(187.8, 75.1)[l]{$e_1$\strut}}
  \put(1413.78,3016.69){\fontsize{10.69}{12.82}\selectfont \makebox(187.8, 75.1)[l]{$e_2$\strut}}
  \put(1202.54,2865.58){\fontsize{10.69}{12.83}\selectfont \makebox(187.8, 75.1)[l]{$G_1$\strut}}
  \put(1569.50,3119.75){\fontsize{10.69}{12.82}\selectfont \makebox(112.7, 75.1)[l]{$z$\strut}}
  \put(2249.04,2758.93){\fontsize{10.69}{12.82}\selectfont \makebox(112.7, 75.1)[l]{(f)\strut}}
  \put(2267.82,3173.56){\fontsize{10.69}{12.82}\selectfont \makebox(187.8, 75.1)[l]{$e_1$\strut}}
  \put(2253.28,2987.93){\fontsize{10.69}{12.82}\selectfont \makebox(187.8, 75.1)[l]{$e_2$\strut}}
  \put(2046.50,2864.95){\fontsize{10.69}{12.82}\selectfont \makebox(187.8, 75.1)[l]{$G_1$\strut}}
  \put(442.54,2101.74){\fontsize{10.69}{12.83}\selectfont \makebox(112.7, 75.1)[l]{(g)\strut}}
  \put(393.20,2574.50){\fontsize{10.69}{12.82}\selectfont \makebox(187.8, 75.1)[l]{$e_1$\strut}}
  \put(390.29,2297.86){\fontsize{10.69}{12.82}\selectfont \makebox(187.8, 75.1)[l]{$e_2$\strut}}
  \put(161.52,2190.55){\fontsize{10.69}{12.82}\selectfont \makebox(187.8, 75.1)[l]{$G_1$\strut}}
  \put(1742.06,3163.74){\fontsize{10.69}{12.82}\selectfont \makebox(112.7, 75.1)[l]{$x$\strut}}
  \put(1742.06,3075.76){\fontsize{10.69}{12.82}\selectfont \makebox(112.7, 75.1)[l]{$y$\strut}}
  \put(2572.00,3127.27){\fontsize{10.69}{12.82}\selectfont \makebox(112.7, 75.1)[l]{$y$\strut}}
  \put(2409.58,3181.41){\fontsize{10.69}{12.82}\selectfont \makebox(112.7, 75.1)[l]{$x$\strut}}
  \put(2565.24,3035.91){\fontsize{10.69}{12.82}\selectfont \makebox(112.7, 75.1)[l]{$z$\strut}}
  \put(707.49,2527.18){\fontsize{10.69}{12.82}\selectfont \makebox(112.7, 75.1)[l]{$y$\strut}}
  \put(714.12,2429.05){\fontsize{10.69}{12.82}\selectfont \makebox(112.7, 75.1)[l]{$z$\strut}}
  \put(517.45,2465.93){\fontsize{10.69}{12.82}\selectfont \makebox(112.7, 75.1)[l]{$x$\strut}}
  \put(172.62,2832.19){\fontsize{10.69}{12.82}\selectfont \makebox(187.8, 75.1)[l]{$G_1$\strut}}
  \put(1108.94,3506.82){\fontsize{10.69}{12.82}\selectfont \makebox(187.8, 75.1)[l]{$G_1$\strut}}
  \put(2029.00,3503.61){\fontsize{10.69}{12.82}\selectfont \makebox(187.8, 75.1)[l]{$G_1$\strut}}
  \put(437.43,3198.09){\fontsize{10.69}{12.82}\selectfont \makebox(187.8, 75.1)[l]{$e_1$\strut}}
  \put(437.83,3014.48){\fontsize{10.69}{12.82}\selectfont \makebox(187.8, 75.1)[l]{$e_2$\strut}}
  \put(2267.20,3889.44){\fontsize{10.69}{12.82}\selectfont \makebox(187.8, 75.1)[l]{$e_1$\strut}}
  \put(2291.46,3679.52){\fontsize{10.69}{12.82}\selectfont \makebox(187.8, 75.1)[l]{$e_2$\strut}}
  \put(2577.18,3819.06){\fontsize{10.69}{12.82}\selectfont \makebox(187.8, 75.1)[l]{$e_3$\strut}}
  \put(1363.78,3882.18){\fontsize{10.69}{12.82}\selectfont \makebox(187.8, 75.1)[l]{$e_1$\strut}}
  \put(1343.85,3676.20){\fontsize{10.69}{12.82}\selectfont \makebox(187.8, 75.1)[l]{$e_2$\strut}}
  \put(1602.99,3785.83){\fontsize{10.69}{12.82}\selectfont \makebox(187.8, 75.1)[l]{$e_3$\strut}}
  \put(714.16,3150.96){\fontsize{10.69}{12.82}\selectfont \makebox(187.8, 75.1)[l]{$e_3$\strut}}
  \put(720.81,2984.85){\fontsize{10.69}{12.82}\selectfont \makebox(187.8, 75.1)[l]{$e_4$\strut}}
  \put(352.54,4007.17){\fontsize{10.69}{12.82}\selectfont \makebox(112.7, 75.1)[l]{$u$\strut}}
  \put(385.62,3406.03){\fontsize{10.69}{12.82}\selectfont \makebox(112.7, 75.1)[l]{$v$\strut}}
  \put(1271.69,2401.43){\fontsize{10.69}{12.82}\selectfont \makebox(112.7, 75.1)[l]{$P$\strut}}
  \put(907.29,2401.63){\fontsize{10.69}{12.82}\selectfont \makebox(187.8, 75.1)[l]{$P_1$\strut}}
  \put(1436.92,2042.13){\fontsize{10.69}{12.83}\selectfont \makebox(112.7, 75.1)[l]{(h)\strut}}
  \put(1221.92,2705.51){\fontsize{10.69}{12.82}\selectfont \makebox(112.7, 75.1)[l]{$u$\strut}}
  \put(1231.69,2118.13){\fontsize{10.69}{12.82}\selectfont \makebox(112.7, 75.1)[l]{$v$\strut}}
  \put(1477.83,2561.94){\fontsize{10.69}{12.83}\selectfont \makebox(187.8, 75.1)[l]{$e_3$\strut}}
  \put(1508.49,2245.11){\fontsize{10.69}{12.83}\selectfont \makebox(187.8, 75.1)[l]{$e_4$\strut}}
  \put(400.13,1929.34){\fontsize{10.69}{12.82}\selectfont \makebox(112.7, 75.1)[l]{$u$\strut}}
  \put(398.94,1521.13){\fontsize{10.69}{12.82}\selectfont \makebox(112.7, 75.1)[l]{$v$\strut}}
  \put(179.97,1733.11){\fontsize{10.69}{12.82}\selectfont \makebox(112.7, 75.1)[l]{$x$\strut}}
  \put(632.27,1809.49){\fontsize{10.69}{12.82}\selectfont \makebox(112.7, 75.1)[l]{$y$\strut}}
  \put(450.08,1712.57){\fontsize{10.69}{12.82}\selectfont \makebox(112.7, 75.1)[l]{$z$\strut}}
  \put(1547.60,1371.79){\fontsize{10.69}{12.83}\selectfont \makebox(112.7, 75.1)[l]{(k)\strut}}
  \put(530.04,1358.02){\fontsize{10.69}{12.83}\selectfont \makebox(112.7, 75.1)[l]{(j)\strut}}
  \put(282.47,1410.13){\fontsize{10.69}{12.82}\selectfont \makebox(187.8, 75.1)[l]{$G_1$\strut}}
  \put(1367.88,1925.54){\fontsize{10.69}{12.82}\selectfont \makebox(112.7, 75.1)[l]{$u$\strut}}
  \put(1386.19,1515.01){\fontsize{10.69}{12.82}\selectfont \makebox(112.7, 75.1)[l]{$v$\strut}}
  \put(1183.62,1731.82){\fontsize{10.69}{12.82}\selectfont \makebox(112.7, 75.1)[l]{$x$\strut}}
  \put(2325.70,1788.45){\fontsize{10.69}{12.82}\selectfont \makebox(112.7, 75.1)[l]{$y$\strut}}
  \put(1444.40,1717.64){\fontsize{10.69}{12.82}\selectfont \makebox(112.7, 75.1)[l]{$z$\strut}}
  \put(1280.44,1396.98){\fontsize{10.69}{12.82}\selectfont \makebox(187.8, 75.1)[l]{$G_1$\strut}}
  \put(2258.56,1427.63){\fontsize{10.69}{12.83}\selectfont \makebox(112.7, 75.1)[l]{(l)\strut}}
  \put(376.58,1253.13){\fontsize{10.69}{12.82}\selectfont \makebox(112.7, 75.1)[l]{$u$\strut}}
  \put(391.81,830.10){\fontsize{10.69}{12.82}\selectfont \makebox(112.7, 75.1)[l]{$v$\strut}}
  \put(180.35,1052.12){\fontsize{10.69}{12.82}\selectfont \makebox(112.7, 75.1)[l]{$x$\strut}}
  \put(612.00,1120.96){\fontsize{10.69}{12.82}\selectfont \makebox(112.7, 75.1)[l]{$y$\strut}}
  \put(444.39,1027.68){\fontsize{10.69}{12.82}\selectfont \makebox(112.7, 75.1)[l]{$z$\strut}}
  \put(551.49,678.41){\fontsize{10.69}{12.83}\selectfont \makebox(112.7, 75.1)[l]{(m)\strut}}
  \put(276.78,725.24){\fontsize{10.69}{12.82}\selectfont \makebox(187.8, 75.1)[l]{$G_1$\strut}}
  \put(1259.88,1237.76){\fontsize{10.69}{12.82}\selectfont \makebox(112.7, 75.1)[l]{$u$\strut}}
  \put(1266.76,818.61){\fontsize{10.69}{12.82}\selectfont \makebox(112.7, 75.1)[l]{$v$\strut}}
  \put(1039.73,1041.53){\fontsize{10.69}{12.82}\selectfont \makebox(112.7, 75.1)[l]{$x$\strut}}
  \put(1492.03,1117.91){\fontsize{10.69}{12.82}\selectfont \makebox(112.7, 75.1)[l]{$y$\strut}}
  \put(2364.30,680.20){\fontsize{10.69}{12.83}\selectfont \makebox(112.7, 75.1)[l]{(o)\strut}}
  \put(1389.80,687.36){\fontsize{10.69}{12.83}\selectfont \makebox(112.7, 75.1)[l]{(n)\strut}}
  \put(1142.22,718.55){\fontsize{10.69}{12.82}\selectfont \makebox(187.8, 75.1)[l]{$G_1$\strut}}
  \put(2184.60,1233.96){\fontsize{10.69}{12.82}\selectfont \makebox(112.7, 75.1)[l]{$u$\strut}}
  \put(2197.92,819.96){\fontsize{10.69}{12.82}\selectfont \makebox(112.7, 75.1)[l]{$v$\strut}}
  \put(2000.34,1040.23){\fontsize{10.69}{12.82}\selectfont \makebox(112.7, 75.1)[l]{$x$\strut}}
  \put(2436.02,1115.69){\fontsize{10.69}{12.82}\selectfont \makebox(112.7, 75.1)[l]{$y$\strut}}
  \put(2097.16,705.40){\fontsize{10.69}{12.82}\selectfont \makebox(187.8, 75.1)[l]{$G_1$\strut}}
  \put(776.21,121.49){\fontsize{10.69}{12.83}\selectfont \makebox(112.7, 75.1)[l]{(p)\strut}}
  \put(1714.02,616.06){\fontsize{10.69}{12.82}\selectfont \makebox(112.7, 75.1)[l]{$u$\strut}}
  \put(1721.02,190.86){\fontsize{10.69}{12.82}\selectfont \makebox(112.7, 75.1)[l]{$v$\strut}}
  \put(1517.78,415.05){\fontsize{10.69}{12.82}\selectfont \makebox(112.7, 75.1)[l]{$x$\strut}}
  \put(1949.44,483.88){\fontsize{10.69}{12.82}\selectfont \makebox(112.7, 75.1)[l]{$y$\strut}}
  \put(1885.34,37.75){\fontsize{10.69}{12.83}\selectfont \makebox(112.7, 75.1)[l]{(q)\strut}}
  \put(1614.21,88.17){\fontsize{10.69}{12.82}\selectfont \makebox(187.8, 75.1)[l]{$G_1$\strut}}
  \put(660.70,561.21){\fontsize{10.69}{12.82}\selectfont \makebox(112.7, 75.1)[l]{$u$\strut}}
  \put(664.78,198.95){\fontsize{10.69}{12.82}\selectfont \makebox(112.7, 75.1)[l]{$v$\strut}}
  \put(481.66,331.26){\fontsize{10.69}{12.82}\selectfont \makebox(112.7, 75.1)[l]{$x$\strut}}
  \put(827.42,448.23){\fontsize{10.69}{12.82}\selectfont \makebox(112.7, 75.1)[l]{$y$\strut}}
  \put(2256.32,2014.89){\fontsize{10.69}{12.83}\selectfont \makebox(112.7, 75.1)[l]{(i)\strut}}
  \put(2220.40,2396.36){\fontsize{10.69}{12.82}\selectfont \makebox(112.7, 75.1)[l]{$P$\strut}}
  \put(1856.00,2396.56){\fontsize{10.69}{12.82}\selectfont \makebox(187.8, 75.1)[l]{$P_1$\strut}}
  \put(2191.74,2108.17){\fontsize{10.69}{12.82}\selectfont \makebox(112.7, 75.1)[l]{$v$\strut}}
  \put(2428.52,2506.37){\fontsize{10.69}{12.83}\selectfont \makebox(187.8, 75.1)[l]{$e_3$\strut}}
  \put(2249.14,2240.90){\fontsize{10.69}{12.83}\selectfont \makebox(187.8, 75.1)[l]{$e_2$\strut}}
  \put(2152.26,2699.89){\fontsize{10.69}{12.82}\selectfont \makebox(112.7, 75.1)[l]{$u$\strut}}
  \put(2155.24,1887.89){\fontsize{10.69}{12.82}\selectfont \makebox(112.7, 75.1)[l]{$u$\strut}}
  \put(1948.62,1720.52){\fontsize{10.69}{12.82}\selectfont \makebox(112.7, 75.1)[l]{$x$\strut}}
  \put(2198.58,1702.14){\fontsize{10.69}{12.82}\selectfont \makebox(112.7, 75.1)[l]{$z$\strut}}
  \put(2160.40,1536.63){\fontsize{10.69}{12.82}\selectfont \makebox(112.7, 75.1)[l]{$v$\strut}}
  \put(1623.19,1803.95){\fontsize{10.69}{12.82}\selectfont \makebox(112.7, 75.1)[l]{$y$\strut}}
  \end{picture}%
  \caption{The graphs in the proof of Theorem \ref{thm2.7}.} \label{fig5}
\end{figure}

\noindent {\bf Case 1.}  At least one of $P_1$, $P_2$ and $P_3$, say
$P_2$ has length not less than $3$.

Let $e_1$ and $e_2$ be the edges on $P_2$ incident with $u$ and $v$,
respectively. Then $G-\{e_1, e_2\}=G_1+G_2$, where $G_1$ is a
unicyclic graph and $G_2$ is a tree of order at least $2$. It
follows from Lemma \ref{lem1.3} that $E(G_1)\geq |V(G_1)|$. If
$G_2\not\cong S_3, S_4, W, S_2, Q$, then we have $E(G_2)> |V(G_2)|$
by Lemma \ref{lem1.1} and Theorem \ref{thm2.5}, and so $E(G)>n$ by
Lemma \ref{lem2.4}. Hence we only need to consider the following
five subcases.

\noindent {\bf Subcase 1.1.} $G_2\cong S_3$. Then $G$ must have the
structure as given in Figure \ref{fig5} (b) or (c). In either case,
$G-\{e_2,e_3\}=G'_1+G'_2$, where $G'_1$ is a unicyclic graph and
$G'_2\cong S_2$. Obviously, $G'_1\not\cong K_{2,2}$. Then
$E(G'_1)>|V(G'_1)|$ by Lemma \ref{lem1.3} and Theorems \ref{thm2.6}.
Since $E(G'_2)=|V(G'_2)|$, we have $E(G)>n$ by Lemma \ref{lem2.4}.

\noindent {\bf Subcase 1.2.} $G_2\cong S_4$. Then $G$ must have the
structure as given in Figure \ref{fig5} (d). Obviously,
$G-\{e_3,e_4\}=G'_1+G'_2$, where $G'_1$ is a unicyclic graph which
is not isomorphic to $K_{2,2}$ and $G'_2\cong S_2$. Similar to the
proof of Subcase 1.1, we have $E(G)>n$.

\noindent {\bf Subcase 1.3.} $G_2\cong W$. Then $G$ must have the
structure as given in Figure \ref{fig5} (e), (f) or (g). Obviously,
$G-\{xy,yz\}=G'_1+G'_2$, where $G'_1$ is a unicyclic graph which is
not isomorphic to $K_{2,2}$ and $G'_2$ is a tree of order $5$ or
$2$. Similarly, we can obtain that $E(G)>n$.

\noindent {\bf Subcase 1.4.} $G_2\cong S_2$. Since $G_1$ is a
unicyclic graph, if $G_1\not\cong K_{2,2}$, then we can similarly
obtain that $E(G)>n$. If $G_1\cong K_{2,2}$, then $n=6$, which is a
contradiction.

\noindent {\bf Subcase 1.5.} $G_2\cong Q$. Then $G$ must have the
structure as given in Figure \ref{fig5} (h) or (i). In the former
case, $G-\{e_3, e_4\}=G'_1+G'_2$, where $G'_2$ is a path of order
$4$ and $G'_1$ is a unicyclic graph which is not isomorphic to
$K_{2,2}$. Similarly, we can obtain that $E(G)>n$. In the latter
case, $G-\{e_2, e_3\}=G'_1+G'_2$, where $G'_2$ is a tree of order
$5$ and $G'_1$ is a unicyclic graph which is not isomorphic to
$K_{2,2}$. Similarly, we can obtain that $E(G)>n$.

\noindent {\bf Case 2.} All the paths $P_1$, $P_2$ and $P_3$ have
length $2$.

We assume that $P_1=uxv$, $P=uzv$ and $P_2=uyv$. Let $F=\{uy,vy\}$,
then $G-F=G_1+G_2$, where $G_1$ is a unicyclic graph and $G_2$ is a
tree. It follows from Lemma \ref{lem1.3} that $E(G_1)\geq |V(G_1)|$.
If $G_2\not\cong S_1, S_3, S_4, W$, then we have $E(G_2)\geq
|V(G_2)|$ by Lemma \ref{lem1.1} and so $E(G)>n$ by Lemma
\ref{lem2.4}. Hence we only need to consider the following four
subcases.

\noindent {\bf Subcase 2.1.} $G_2\cong S_1$. Let $F'=\{uy,zv,xv\}$,
then $G-F'=G'_1+G'_2$, where $G'_2\cong S_2$ and $G'_1$ is a tree of
order at least $6$ since $n\geq 8$. It is easy to see that $G'_1$
can not be isomorphic to $Q$ or $W$. Therefore we have
$E(G'_1)>|V(G'_1)|$ and $E(G'_2)=|V(G'_2)|$ by Lemma \ref{lem1.1}
and Theorem \ref{thm2.5}. It follows from Lemma \ref{lem2.4} that
$E(G)>n$.

\noindent {\bf Subcase 2.2.} $G_2\cong S_3$. Then $G$ must have the
structure as given in Figure \ref{fig5} (j). Let $F'=\{uy,zv,xv\}$,
then $G-F'=G'_1+G'_2$, where $G'_2$ is the path of order $4$
containing $y$ and $G'_1$ is a tree of order at least $4$ since
$n\geq 8$. Clearly, $G'_1$ can not be isomorphic to $S_4$, $Q$ or
$W$. Similar to the proof of Subcase 2.1, we have $E(G)>n$.

\noindent {\bf Subcase 2.3.} $G_2\cong S_4$. Then $G$ must have the
structure as given in Figure \ref{fig5} (k). Let $F'=\{uy,zv,xv\}$,
then $G-F'=G'_1+G'_2$, where $G'_2$ is the tree of order $5$
containing $y$ and $G'_1$ is a tree of order at least $3$. Clearly,
$G'_1$ can not be isomorphic to $S_4$ or $W$. If $G'_1\not\cong
S_3$, then we can similarly obtain that $E(G)>n$. If $G'_1\cong
S_3$, then $G$ must be the graph as given in Figure \ref{fig5} (l).
By choosing the edge cut $\{uy,uz,xv\}$, we can also obtain that
$E(G)>n$.

\noindent {\bf Subcase 2.4.} $G_2\cong W$. Then $G$ must have the
structure as given in Figure \ref{fig5} (m). Let $F'=\{uy,zv,xv\}$,
then $G-F'=G'_1+G'_2$, where $G'_2$ is the tree of order $8$
containing $y$ and $G'_1$ is a tree of order at least $3$. Clearly,
$G'_1$ can not be isomorphic to $S_4$ or $W$. If $G'_1\cong S_3$,
then $n=11$, which is a contradiction. If $G'_1\not\cong S_3$, then
we can similarly obtain that $E(G)>n$.

\noindent {\bf Case 3.} One of the paths $P_1$, $P_2$ and $P_3$ has
length $1$, and the other two paths have length $2$.

Without loss of generality, we assume that $P=uv$, $P_1=uxv$ and
$P_2=uyv$. Let $F=\{uy,vy\}$, then $G-F=G_1+G_2$, where $G_1$ is a
unicyclic graph and $G_2$ is a tree. Similarly, if $G_2\not\cong
S_1, S_3, S_4, W$, then we have $E(G)>n$. Hence we also need to
consider the following four subcases.

\noindent {\bf Subcase 3.1.} $G_2\cong S_1$. Let $F'=\{uy,uv,xv\}$,
then $G-F'=G'_1+G'_2$, where $G'_2\cong S_2$ and $G'_1$ is a tree of
order at least $6$ since $n\geq 8$. Since $\Delta(G)\leq 3$, $G'_1$
can not be isomorphic to $Q$ or $W$. Similar to the proof of Subcase
2.1, we have $E(G)>n$.

\noindent {\bf Subcase 3.2.} $G_2\cong S_3$. Then $G$ must have the
structure as given in Figure \ref{fig5} (n). Let $F'=\{uy,uv,xv\}$,
then $G-F'=G'_1+G'_2$, where $G'_2$ is the path of order $4$
containing $y$ and $G'_1$ is a tree of order at least $4$ since
$n\geq 8$. Clearly, $G'_1$ can not be isomorphic to $S_4$ or $W$.
Similarly, we have $E(G)>n$.

\noindent {\bf Subcase 3.3.} $G_2\cong S_4$. Then $G$ must have the
structure as given in Figure \ref{fig5} (o). Let $F'=\{uy,uv,xv\}$,
then $G-F'=G'_1+G'_2$, where $G'_2$ is the tree of order $5$
containing $y$ and $G'_1$ is a tree of order at least $3$. Clearly,
$G'_1$ can not be isomorphic to $S_4$ or $W$. If $G'_1\not\cong
S_3$, then we can similarly obtain that $E(G)>n$. If $G'_1\cong
S_3$, then $G$ must be the graph as given in Figure \ref{fig5} (p).
By choosing the edge cut $\{xu,xv\}$, we can similarly obtain that
$E(G)>n$.

\noindent {\bf Subcase 3.4.} $G_2\cong W$. Then $G$ must have the
structure as given in Figure \ref{fig5} (q). Let $F'=\{uy,uv,xv\}$,
then $G-F'=G'_1+G'_2$, where $G'_2$ is the tree of order $8$
containing $y$ and $G'_1$ is a tree of order at least $2$. Clearly,
$G'_1$ can not be isomorphic to $S_4$ or $W$. If $G'_1\cong S_3$,
then $n=11$, which is a contradiction. If $G'_1\not\cong S_3$, then
we can similarly obtain that $E(G)>n$. The proof is thus complete.
\qed \\[3mm]
\noindent{\bf Proof of Conjecture \ref{con1.5}:} Let $G$ be a
connected graph of order $n$ with $\Delta \leq 3$. Clearly, if $G$
is isomorphic to a graph in $\{S_2, Q, K_{2,2},K_{3,3}\}$, then
$E(G)=n$. We will prove that $E(G)\neq n$ if $G\not\cong S_2$, $Q$,
$K_{2,2}$ or $K_{3,3}$ by induction on the cyclomatic number $c(G)$.
It follows from Theorems \ref{thm2.5}, \ref{thm2.6} and \ref{thm2.7}
that the result holds for $c(G)\leq 2$. Let $k\geq 3$ be an integer.
We assume that the result holds for $c(G)<k$. Now let $G$ be a graph
with $c(G)=k\geq 3$. We will show that $E(G)\neq n$.

By Lemma \ref{lem2.3}, the result holds if $n$ is odd. By the fact
that $K_{3,3}$ is the only connected cycle-containing graph of order
6 with $\Delta \leq 3$ and $E=6$, we know that the result holds for
$n\leq 6$. So in the following we assume that $n\geq 8$ is even. In
our proof we will repeatedly make use of the following claim:

\noindent {\bf Claim 1.} {\it Let $F$ be an edge cut of $G$ such
that $G-F=G_1 + G_2$ with $c(G_1), c(G_2)<k$. If $G_1, G_2\not\cong
S_1, S_3, S_4, W$ or $K_{2,3}$ and either the edges in $F$ form a
star or at least one of $G_1$ and $G_2$ is not isomorphic to $S_2,
Q$ or $K_{2,2}$, then we are done. }

\pf By Lemma \ref{lem1.4}, we have $E(G_1)\geq |V(G_1)|$ and
$E(G_2)\geq |V(G_2)|$. Clearly, $G_1, G_2\not\cong K_{3,3}$. If
$G_i\not\cong S_2, Q$ or $K_{2,2}$, then by induction hypothesis, we
have $E(G_i)\neq |V(G_i)|$. Therefore we have $E(G)>n$ by Lemma
\ref{lem2.4}. \qed

In what follows, we use $\hat{G}$ to denote the graph obtained from
$G$ by repeatedly deleting the pendent vertices. Clearly,
$c(\hat{G})=c(G)$. Denote by $\kappa'(\hat{G})$ the edge
connectivity of $\hat{G}$. Since $\Delta(\hat{G})\leq 3$, we have
$1\leq \kappa'(\hat{G})\leq 3$. Therefore we only need to consider
the following three cases.

\noindent {\bf Case 1.} $\kappa'(\hat{G})=1$.

\begin{figure}[ht]
\centering
   \setlength{\unitlength}{0.05 mm}%
  \begin{picture}(2240.7, 1075.8)(0,0)
  \put(0,0){\includegraphics{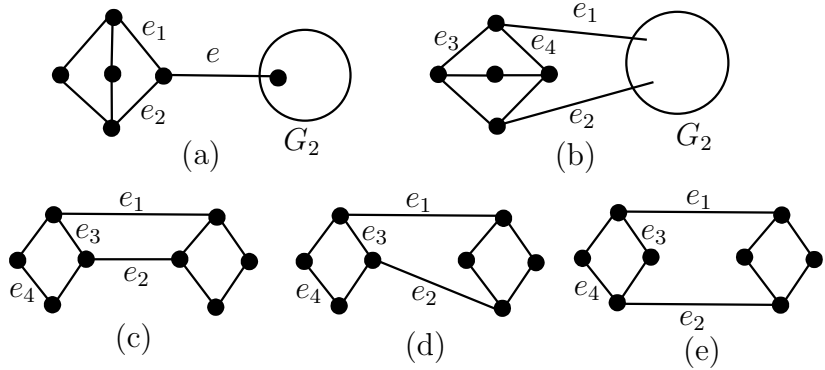}}
  \put(562.97,831.73){\fontsize{11.76}{14.11}\selectfont \makebox(124.0, 82.6)[l]{$e$\strut}}
  \put(392.02,902.47){\fontsize{11.76}{14.11}\selectfont \makebox(206.6, 82.6)[l]{$e_1$\strut}}
  \put(389.07,678.47){\fontsize{11.76}{14.11}\selectfont \makebox(206.6, 82.6)[l]{$e_2$\strut}}
  \put(770.16,596.85){\fontsize{11.76}{14.11}\selectfont \makebox(206.6, 82.6)[l]{$G_2$\strut}}
  \put(1168.43,876.79){\fontsize{11.76}{14.11}\selectfont \makebox(206.6, 82.6)[l]{$e_3$\strut}}
  \put(1538.69,950.67){\fontsize{11.76}{14.11}\selectfont \makebox(206.6, 82.6)[l]{$e_1$\strut}}
  \put(1528.28,669.79){\fontsize{11.76}{14.11}\selectfont \makebox(206.6, 82.6)[l]{$e_2$\strut}}
  \put(1813.06,609.31){\fontsize{11.76}{14.11}\selectfont \makebox(206.6, 82.6)[l]{$G_2$\strut}}
  \put(1426.50,865.19){\fontsize{11.76}{14.11}\selectfont \makebox(206.6, 82.6)[l]{$e_4$\strut}}
  \put(498.20,559.41){\fontsize{11.76}{14.11}\selectfont \makebox(124.0, 82.6)[l]{(a)\strut}}
  \put(1491.20,556.93){\fontsize{11.76}{14.11}\selectfont \makebox(124.0, 82.6)[l]{(b)\strut}}
  \put(332.21,447.40){\fontsize{11.76}{14.11}\selectfont \makebox(206.6, 82.6)[l]{$e_1$\strut}}
  \put(343.14,248.57){\fontsize{11.76}{14.11}\selectfont \makebox(206.6, 82.6)[l]{$e_2$\strut}}
  \put(216.17,362.62){\fontsize{11.76}{14.11}\selectfont \makebox(206.6, 82.6)[l]{$e_3$\strut}}
  \put(40.00,205.28){\fontsize{11.76}{14.11}\selectfont \makebox(206.6, 82.6)[l]{$e_4$\strut}}
  \put(1094.30,444.50){\fontsize{11.76}{14.11}\selectfont \makebox(206.6, 82.6)[l]{$e_1$\strut}}
  \put(1109.38,197.32){\fontsize{11.76}{14.11}\selectfont \makebox(206.6, 82.6)[l]{$e_2$\strut}}
  \put(977.61,358.38){\fontsize{11.76}{14.11}\selectfont \makebox(206.6, 82.6)[l]{$e_3$\strut}}
  \put(802.09,202.38){\fontsize{11.76}{14.11}\selectfont \makebox(206.6, 82.6)[l]{$e_4$\strut}}
  \put(1839.02,453.55){\fontsize{11.76}{14.11}\selectfont \makebox(206.6, 82.6)[l]{$e_1$\strut}}
  \put(1821.93,131.44){\fontsize{11.76}{14.11}\selectfont \makebox(206.6, 82.6)[l]{$e_2$\strut}}
  \put(1719.44,360.82){\fontsize{11.76}{14.11}\selectfont \makebox(206.6, 82.6)[l]{$e_3$\strut}}
  \put(1541.83,210.79){\fontsize{11.76}{14.11}\selectfont \makebox(206.6, 82.6)[l]{$e_4$\strut}}
  \put(321.18,72.87){\fontsize{11.76}{14.11}\selectfont \makebox(124.0, 82.6)[l]{(c)\strut}}
  \put(1092.30,50.91){\fontsize{11.76}{14.11}\selectfont \makebox(124.0, 82.6)[l]{(d)\strut}}
  \put(1832.73,37.52){\fontsize{11.76}{14.11}\selectfont \makebox(124.0, 82.6)[l]{(e)\strut}}
  \end{picture}%
  \caption{The graphs in the proof of Case 1 and Subcase 2.1 of Conjecture \ref{con1.5}.} \label{fig6}
\end{figure}

Let $e$ be a cut edge of $\hat{G}$. Then $\hat{G}-e$ has exactly two
components, say, $H_1$ and $H_2$. It is clear that $c(H_1)\geq 1$,
$c(H_2)\geq 1$ and $c(H_1)+c(H_2)=k$. Consequently, $G-e$ has
exactly two components $G_1$ and $G_2$ with $c(G_1)\geq 1$,
$c(G_2)\geq 1$ and $c(G_1)+c(G_2)=k$, where $H_i$ is a subgraph of
$G_i$ for $i=1,2$. If neither $G_1$ nor $G_2$ is isomorphic to
$K_{2,3}$, then we are done by Claim 1. Otherwise, without loss of
generality, we assume that $G_1\cong K_{2,3}$. Then $G$ must have
the structure as given in Figure \ref{fig6} (a). Now, let
$F=\{e_1,e_2\}$. Then $G-F=G'_1+G'_2$, where $G'_1\cong K_{2,2}$ and
$G'_2=G_2\cup e$. Therefore we have that $c(G'_2)=k-2\geq 1$ and
$G'_2\not\cong K_{2,2}, K_{2,3}$, and so we are done by Claim 1.

\noindent {\bf Case 2.} $\kappa'(\hat{G})=2$.

Let $F=\{e_1, e_2\}$ be an edge cut of $\hat{G}$. Then $\hat{G}-F$
has exactly two components, say, $H_1$ and $H_2$. Clearly,
$c(H_1)+c(H_2)=k-1\geq 2$.

\noindent {\bf Subcase 2.1.} $c(H_1)\geq 1$ and $c(H_2)\geq 1$.
Therefore, $G-F$ has exactly two components $G_1$ and $G_2$ with
$c(G_1)\geq 1$, $c(G_2)\geq 1$ and $c(G_1)+c(G_2)=k-1$, where $H_i$
is a subgraph of $G_i$ for $i=1,2$. If $G_1, G_2\not\cong K_{2,3}$
and at least one of $G_1$ and $G_2$ is not isomorphic to $K_{2,2}$,
then we are done by Claim 1. If at least one of $G_1$ and $G_2$ is
isomorphic to $K_{2,3}$, say $G_1\cong K_{2,3}$. Then $G$ must have
the structure as given in Figure \ref{fig6} (b). Now, let
$F'=\{e_2,e_3,e_4\}$, then $G-F'=G'_1+G'_2$, where $G'_1\cong
K_{2,2}$ and $G'_2=G_2\cup e_1$. Therefore we have that
$c(G'_2)=k-3$ and $G'_2\not\cong K_{2,2}, K_{2,3}$, and so we are
done by Claim 1. If $G_1,G_2\cong K_{2,2}$, then $G$ must be the
graph as given in Figure \ref{fig6} (c), (d) or (e). Let $F'=\{e_1,
e_3, e_4\}$, then $G-F'=G'_1+G'_2$, where $G'_1\cong S_2$ and
$G'_2\not\cong K_{2,2}$ is a unicyclic graph. Hence we are done by
Claim 1.

\noindent {\bf Subcase 2.2.} One of $H_1$ and $H_2$, say $H_2$ is a
tree. Therefore, $G-F$ has exactly two components $G_1$ and $G_2$
with $c(G_1)=k-1$ and $c(G_2)=0$, where $H_i$ is a subgraph of $G_i$
for $i=1,2$. Since $k-1\geq 2$, $G_1\not\cong S_2,Q,K_{2,2}$. If
$G_1\not\cong K_{2,3}$ and $G_2\not\cong S_1,S_3,S_4,W$, then we are
done by Claim 1. So we assume that this is not true. We only need to
consider the following five subsubcases.

\begin{figure}[ht]
\centering
    \setlength{\unitlength}{0.05 mm}%
  \begin{picture}(2348.1, 3013.8)(0,0)
  \put(0,0){\includegraphics{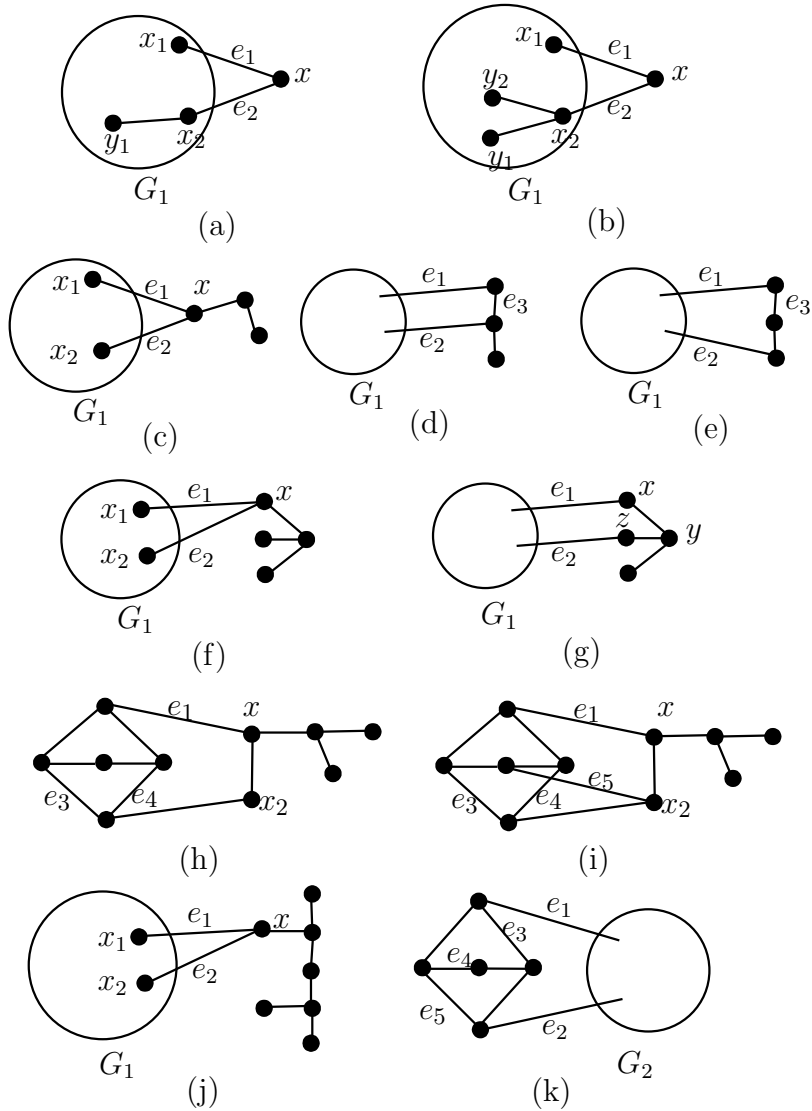}}
  \put(796.17,2749.79){\fontsize{11.76}{14.11}\selectfont \makebox(124.0, 82.6)[l]{$x$\strut}}
  \put(383.52,2841.44){\fontsize{11.76}{14.11}\selectfont \makebox(206.6, 82.6)[l]{$x_1$\strut}}
  \put(483.15,2586.51){\fontsize{11.76}{14.11}\selectfont \makebox(206.6, 82.6)[l]{$x_2$\strut}}
  \put(627.38,2810.25){\fontsize{11.76}{14.11}\selectfont \makebox(206.6, 82.6)[l]{$e_1$\strut}}
  \put(632.41,2661.61){\fontsize{11.76}{14.11}\selectfont \makebox(206.6, 82.6)[l]{$e_2$\strut}}
  \put(1799.89,2755.59){\fontsize{11.76}{14.11}\selectfont \makebox(124.0, 82.6)[l]{$x$\strut}}
  \put(1383.95,2859.51){\fontsize{11.76}{14.11}\selectfont \makebox(206.6, 82.6)[l]{$x_1$\strut}}
  \put(1478.60,2586.98){\fontsize{11.76}{14.11}\selectfont \makebox(206.6, 82.6)[l]{$x_2$\strut}}
  \put(1628.58,2823.61){\fontsize{11.76}{14.11}\selectfont \makebox(206.6, 82.6)[l]{$e_1$\strut}}
  \put(1626.06,2664.89){\fontsize{11.76}{14.11}\selectfont \makebox(206.6, 82.6)[l]{$e_2$\strut}}
  \put(538.57,2343.69){\fontsize{11.76}{14.11}\selectfont \makebox(124.0, 82.6)[l]{(a)\strut}}
  \put(1572.37,2352.81){\fontsize{11.76}{14.11}\selectfont \makebox(124.0, 82.6)[l]{(b)\strut}}
  \put(369.03,2437.55){\fontsize{11.76}{14.11}\selectfont \makebox(206.6, 82.6)[l]{$G_1$\strut}}
  \put(1363.67,2437.46){\fontsize{11.76}{14.11}\selectfont \makebox(206.6, 82.6)[l]{$G_1$\strut}}
  \put(290.36,2573.54){\fontsize{11.76}{14.11}\selectfont \makebox(206.6, 82.6)[l]{$y_1$\strut}}
  \put(1311.14,2528.85){\fontsize{11.76}{14.11}\selectfont \makebox(206.6, 82.6)[l]{$y_1$\strut}}
  \put(1296.02,2748.03){\fontsize{11.76}{14.11}\selectfont \makebox(206.6, 82.6)[l]{$y_2$\strut}}
  \put(528.29,2186.47){\fontsize{11.76}{14.11}\selectfont \makebox(124.0, 82.6)[l]{$x$\strut}}
  \put(153.15,2197.19){\fontsize{11.76}{14.11}\selectfont \makebox(206.6, 82.6)[l]{$x_1$\strut}}
  \put(145.11,2012.30){\fontsize{11.76}{14.11}\selectfont \makebox(206.6, 82.6)[l]{$x_2$\strut}}
  \put(396.92,2186.32){\fontsize{11.76}{14.11}\selectfont \makebox(206.6, 82.6)[l]{$e_1$\strut}}
  \put(401.96,2037.68){\fontsize{11.76}{14.11}\selectfont \makebox(206.6, 82.6)[l]{$e_2$\strut}}
  \put(389.68,1781.31){\fontsize{11.76}{14.11}\selectfont \makebox(124.0, 82.6)[l]{(c)\strut}}
  \put(206.74,1848.85){\fontsize{11.76}{14.11}\selectfont \makebox(206.6, 82.6)[l]{$G_1$\strut}}
  \put(1134.52,2218.72){\fontsize{11.76}{14.11}\selectfont \makebox(206.6, 82.6)[l]{$e_1$\strut}}
  \put(1128.15,2044.30){\fontsize{11.76}{14.11}\selectfont \makebox(206.6, 82.6)[l]{$e_2$\strut}}
  \put(1101.08,1814.47){\fontsize{11.76}{14.11}\selectfont \makebox(124.0, 82.6)[l]{(d)\strut}}
  \put(941.55,1893.32){\fontsize{11.76}{14.11}\selectfont \makebox(206.6, 82.6)[l]{$G_1$\strut}}
  \put(1870.34,2219.86){\fontsize{11.76}{14.11}\selectfont \makebox(206.6, 82.6)[l]{$e_1$\strut}}
  \put(1856.84,2006.69){\fontsize{11.76}{14.11}\selectfont \makebox(206.6, 82.6)[l]{$e_2$\strut}}
  \put(1853.68,1802.31){\fontsize{11.76}{14.11}\selectfont \makebox(124.0, 82.6)[l]{(e)\strut}}
  \put(1681.67,1896.07){\fontsize{11.76}{14.11}\selectfont \makebox(206.6, 82.6)[l]{$G_1$\strut}}
  \put(1350.53,2147.32){\fontsize{11.76}{14.11}\selectfont \makebox(206.6, 82.6)[l]{$e_3$\strut}}
  \put(2101.52,2145.20){\fontsize{11.76}{14.11}\selectfont \makebox(206.6, 82.6)[l]{$e_3$\strut}}
  \put(279.84,1588.42){\fontsize{11.76}{14.11}\selectfont \makebox(206.6, 82.6)[l]{$x_1$\strut}}
  \put(279.84,1470.51){\fontsize{11.76}{14.11}\selectfont \makebox(206.6, 82.6)[l]{$x_2$\strut}}
  \put(522.01,1197.37){\fontsize{11.76}{14.11}\selectfont \makebox(124.0, 82.6)[l]{(f)\strut}}
  \put(746.08,1636.65){\fontsize{11.76}{14.11}\selectfont \makebox(124.0, 82.6)[l]{$x$\strut}}
  \put(1548.54,656.42){\fontsize{11.76}{14.11}\selectfont \makebox(124.0, 82.6)[l]{(i)\strut}}
  \put(488.53,656.42){\fontsize{11.76}{14.11}\selectfont \makebox(124.0, 82.6)[l]{(h)\strut}}
  \put(1476.39,1647.40){\fontsize{11.76}{14.11}\selectfont \makebox(206.6, 82.6)[l]{$e_1$\strut}}
  \put(1479.09,1474.71){\fontsize{11.76}{14.11}\selectfont \makebox(206.6, 82.6)[l]{$e_2$\strut}}
  \put(1507.28,1208.28){\fontsize{11.76}{14.11}\selectfont \makebox(124.0, 82.6)[l]{(g)\strut}}
  \put(1292.71,1307.65){\fontsize{11.76}{14.11}\selectfont \makebox(206.6, 82.6)[l]{$G_1$\strut}}
  \put(511.75,1644.13){\fontsize{11.76}{14.11}\selectfont \makebox(206.6, 82.6)[l]{$e_1$\strut}}
  \put(514.45,1471.45){\fontsize{11.76}{14.11}\selectfont \makebox(206.6, 82.6)[l]{$e_2$\strut}}
  \put(325.39,1296.34){\fontsize{11.76}{14.11}\selectfont \makebox(206.6, 82.6)[l]{$G_1$\strut}}
  \put(1710.94,1650.63){\fontsize{11.76}{14.11}\selectfont \makebox(124.0, 82.6)[l]{$x$\strut}}
  \put(1836.66,1535.41){\fontsize{11.76}{14.11}\selectfont \makebox(124.0, 82.6)[l]{$y$\strut}}
  \put(1649.10,1564.89){\fontsize{11.76}{14.11}\selectfont \makebox(124.0, 82.6)[l]{$z$\strut}}
  \put(658.88,1060.99){\fontsize{12.93}{15.52}\selectfont \makebox(136.4, 90.9)[l]{$x$\strut}}
  \put(701.15,811.18){\fontsize{12.93}{15.52}\selectfont \makebox(227.3, 90.9)[l]{$x_2$\strut}}
  \put(1762.01,1059.93){\fontsize{12.93}{15.52}\selectfont \makebox(136.4, 90.9)[l]{$x$\strut}}
  \put(1770.83,803.41){\fontsize{12.93}{15.52}\selectfont \makebox(227.3, 90.9)[l]{$x_2$\strut}}
  \put(458.73,1062.32){\fontsize{12.93}{15.52}\selectfont \makebox(227.3, 90.9)[l]{$e_1$\strut}}
  \put(1537.79,1056.73){\fontsize{12.93}{15.52}\selectfont \makebox(227.3, 90.9)[l]{$e_1$\strut}}
  \put(129.14,823.84){\fontsize{12.93}{15.52}\selectfont \makebox(227.3, 90.9)[l]{$e_3$\strut}}
  \put(1214.21,816.38){\fontsize{12.93}{15.52}\selectfont \makebox(227.3, 90.9)[l]{$e_3$\strut}}
  \put(362.38,830.67){\fontsize{12.93}{15.52}\selectfont \makebox(227.3, 90.9)[l]{$e_4$\strut}}
  \put(1436.74,826.96){\fontsize{12.93}{15.52}\selectfont \makebox(227.3, 90.9)[l]{$e_4$\strut}}
  \put(1575.95,875.33){\fontsize{12.93}{15.52}\selectfont \makebox(227.3, 90.9)[l]{$e_5$\strut}}
  \put(273.76,451.08){\fontsize{11.76}{14.11}\selectfont \makebox(206.6, 82.6)[l]{$x_1$\strut}}
  \put(273.76,333.18){\fontsize{11.76}{14.11}\selectfont \makebox(206.6, 82.6)[l]{$x_2$\strut}}
  \put(508.16,41.22){\fontsize{11.76}{14.11}\selectfont \makebox(124.0, 82.6)[l]{(j)\strut}}
  \put(740.00,499.31){\fontsize{11.76}{14.11}\selectfont \makebox(124.0, 82.6)[l]{$x$\strut}}
  \put(511.11,515.15){\fontsize{11.76}{14.11}\selectfont \makebox(206.6, 82.6)[l]{$e_1$\strut}}
  \put(523.92,364.61){\fontsize{11.76}{14.11}\selectfont \makebox(206.6, 82.6)[l]{$e_2$\strut}}
  \put(277.28,111.56){\fontsize{11.76}{14.11}\selectfont \makebox(206.6, 82.6)[l]{$G_1$\strut}}
  \put(1426.43,37.52){\fontsize{11.76}{14.11}\selectfont \makebox(124.0, 82.6)[l]{(k)\strut}}
  \put(1465.97,541.44){\fontsize{11.76}{14.11}\selectfont \makebox(206.6, 82.6)[l]{$e_1$\strut}}
  \put(1454.46,214.94){\fontsize{11.76}{14.11}\selectfont \makebox(206.6, 82.6)[l]{$e_2$\strut}}
  \put(1654.49,112.52){\fontsize{11.76}{14.11}\selectfont \makebox(206.6, 82.6)[l]{$G_2$\strut}}
  \put(1347.29,478.16){\fontsize{11.76}{14.11}\selectfont \makebox(206.6, 82.6)[l]{$e_3$\strut}}
  \put(1201.42,408.47){\fontsize{11.76}{14.11}\selectfont \makebox(206.6, 82.6)[l]{$e_4$\strut}}
  \put(1128.45,255.04){\fontsize{12.93}{15.52}\selectfont \makebox(227.3, 90.9)[l]{$e_5$\strut}}
  \end{picture}%
  \caption{The graphs in the proof of Subcase 2.2 of Conjecture \ref{con1.5}.} \label{fig7}
\end{figure}

\noindent {\bf Subsubcase 2.2.1.} $G_2\cong S_1$. Let
$V(G_2)=\{x\}$, $e_1=xx_1$ and $e_2=xx_2$. It is clear that
$d_{G_1}(x_2)=1$ or $2$. If $d_{G_1}(x_2)=1$, let
$N_{G_1}(x_2)=\{y_1\}$ (see Figure \ref{fig7} (a), where $y_1$ may
be equal to $x_1$). Let $F'=\{e_1, x_2y_1\}$. Then $G-F'=G'_1+G'_2$,
where $G'_1$ is a graph obtained from $G_1$ by deleting a pendent
vertex and $G'_2\cong S_2$. Therefore, $c(G'_1)=k-1\geq 2$. If
$G'_1\not\cong K_{2,3}$, then we are done by Claim 1. Otherwise,
$n=7$, which is a contradiction.

If $d_{G_1}(x_2)=2$, let $N_{G_1}(x_2)=\{y_1,y_2\}$ (see Figure
\ref{fig7} (b), where one of $y_1$ and $y_2$ may be equal to $x_1$).
Let $F'=\{e_1, x_2y_1, x_2y_2\}$. Then $G-F'=G'_1+G'_2$, where
$G'_1$ is a graph obtained from $G_1$ by deleting a vertex of degree
$2$ and $G'_2\cong S_2$. Therefore, $c(G'_1)=k-2\geq 1$. If
$G'_1\not\cong K_{2,2}, K_{2,3}$, then we are done by Claim 1.
Otherwise, $n=6$ or $7$, which is a contradiction.

\noindent {\bf Subsubcase 2.2.2.} $G_2\cong S_3$. If $e_1$, $e_2$
are incident with a common vertex in $G_2$, then $G$ must have the
structure as given in Figure \ref{fig7} (c). Similar to the proof of
Subsubcase 2.2.1, we can obtain that there exists an edge cut $F'$
such that $G-F'=G'_1+G'_2$ satisfying that $c(G'_1)=k-1$ if
$d_{G_1}(x_2)=1$ or $c(G'_1)=k-2$ if $d_{G_1}(x_2)=2$ and $G'_2$ is
a path of order $4$. If $G'_1\not\cong K_{2,3}$, then we are done by
Claim 1. Otherwise $n=9$, which is a contradiction.

If $e_1$, $e_2$ are incident with two different vertices in $G_2$,
then $G$ must have the structure as given in Figure \ref{fig7} (d)
or (e). Let $F'=\{e_2, e_3\}$, then $G-F'=G'_1+G'_2$, where
$G'_1=G_1\cup e_1$ and $G'_2\cong S_2$. Therefore we have that
$c(G'_1)=k-1\geq 2$ and $G'_2\not\cong K_{2,3}$, and so we are done
by Claim 1.

\noindent {\bf Subsubcase 2.2.3.} $G_2\cong S_4$. If $e_1$, $e_2$
are incident with a common vertex in $G_2$, then $G$ must have the
structure as given in Figure \ref{fig7} (f). Similar to the proof of
Subsubcase 2.2.1, we can obtain that there exists an edge cut $F'$
such that $G-F'=G'_1+G'_2$ satisfying that $c(G'_1)=k-1$ if
$d_{G_1}(x_2)=1$ or $c(G'_1)=k-2$ if $d_{G_1}(x_2)=2$ and $G'_2$ is
a tree of order $5$. If $G'_1\not\cong K_{2,3}$, then we are done by
Claim 1. Otherwise $G$ must be the graph as given in Figure
\ref{fig7} (h) or (i). In the former case let $F''=\{e_1,e_3,e_4\}$
while in the latter case let $F''=\{e_1,e_3,e_4,e_5\}$. Then
$G-F''=G''_1+G''_2$, where $G''_1\cong K_{2,2}$, $G''_2$ is a tree
of order $6$ and $G''_2\not\cong Q$. Therefore we are done by Claim
1.

If $e_1$, $e_2$ are incident with two different vertices in $G_2$,
then $G$ must have the structure as given in Figure \ref{fig7} (g).
Let $F'=\{xy, yz\}$, then $G-F'=G'_1+G'_2$, where $G'_1=G_1\cup
\{e_1,e_2\}$ and $G'_2\cong S_2$. Therefore we have that
$c(G'_1)=k-1\geq 2$ and $G'_2\not\cong K_{2,3}$, and so we are done
by Claim 1.

\noindent {\bf Subsubcase 2.2.4.} $G_2\cong W$. If $e_1$, $e_2$ are
incident with a common vertex in $G_2$, then $G$ must have the
structure as given in Figure \ref{fig7} (j). Similar to the proof of
Subsubcase 2.2.1, we can obtain that there exists an edge cut $F'$
such that $G-F'=G'_1+G'_2$ satisfying that $c(G'_1)=k-1$ if
$d_{G_1}(x_2)=1$ or $c(G'_1)=k-2$ if $d_{G_1}(x_2)=2$ and $G'_2$ is
a tree of order $8$. If $G'_1\not\cong K_{2,3}$, then we are done by
Claim 1. Otherwise, $n=13$, which is a contradiction.

If $e_1$, $e_2$ are incident with two different vertices in $G_2$,
then $G$ must have the structure as given in Figure \ref{fig5} (e),
(f) or (g) ($e_1$, $e_2$ may be incident with a common vertex in
$G_1$). Let $F'=\{xy, yz\}$, then $G-F'=G'_1+G'_2$, where $G'_2$ is
the tree of order $5$ or $2$ containing $y$. Clearly,
$c(G'_1)=k-1\geq 2$ and $G'_1\not\cong K_{2,3}$. Therefore we are
done by Claim 1.

\noindent {\bf Subsubcase 2.2.5.} $G_1\cong K_{2,3}$ and
$G_2\not\cong S_1, S_3, S_4, W$. It is easy to see that $G$ must
have the structure as given in Figure \ref{fig7} (k). Let
$F'=\{e_1,e_3,e_4,e_5\}$. Then $G-F'=G'_1+G'_2$, where $G'_1\cong
S_2$ and $G'_2$ is a tree of order at least $6$ since $n\geq 8$. It
is easy to see that $G'_2$ can not be isomorphic to $W$ or $Q$.
Therefore we are done by Claim 1.

\noindent {\bf Case 3.} $\kappa'(\hat{G})=3$.

Noticing that $\Delta (\hat G)\leq 3$ and $\Delta (G)\leq 3$, we
obtain that $G=\hat G$ is a connected $3$-regular graph. Hence we
have $n+k-1=m=\frac{3}{2}n$, i.e., $n=2k-2$. Since $n\geq 8$, we
have $k\geq 5$.

Let $F=\{e_1, e_2, e_3\}$ be an edge cut of $G$. Then $G-F$ has
exactly two components, say, $G_1$ and $G_2$. Clearly,
$c(G_1)+c(G_2)=k-2\geq 3$. Let $c(G_1)\geq c(G_2)$. If $c(G_2)\geq
3$, then we are done by Claim 1. Hence we only need to consider the
following three subcases.

\begin{figure}[ht]
\centering
  \setlength{\unitlength}{0.05 mm}%
  \begin{picture}(2554.8, 753.5)(0,0)
  \put(0,0){\includegraphics{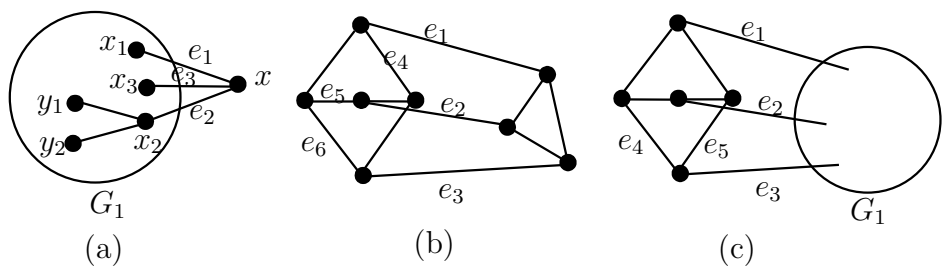}}
  \put(1884.03,316.43){\fontsize{11.76}{14.11}\selectfont \makebox(206.6, 82.6)[l]{$e_5$\strut}}
  \put(1982.75,628.36){\fontsize{11.76}{14.11}\selectfont \makebox(206.6, 82.6)[l]{$e_1$\strut}}
  \put(2021.27,202.48){\fontsize{11.76}{14.11}\selectfont \makebox(206.6, 82.6)[l]{$e_3$\strut}}
  \put(1656.51,329.03){\fontsize{11.76}{14.11}\selectfont \makebox(206.6, 82.6)[l]{$e_4$\strut}}
  \put(2027.59,433.12){\fontsize{11.76}{14.11}\selectfont \makebox(206.6, 82.6)[l]{$e_2$\strut}}
  \put(2273.70,140.38){\fontsize{11.76}{14.11}\selectfont \makebox(206.6, 82.6)[l]{$G_1$\strut}}
  \put(689.72,498.36){\fontsize{11.76}{14.11}\selectfont \makebox(124.0, 82.6)[l]{$x$\strut}}
  \put(277.16,591.64){\fontsize{11.76}{14.11}\selectfont \makebox(206.6, 82.6)[l]{$x_1$\strut}}
  \put(368.43,329.75){\fontsize{11.76}{14.11}\selectfont \makebox(206.6, 82.6)[l]{$x_2$\strut}}
  \put(512.61,567.00){\fontsize{11.76}{14.11}\selectfont \makebox(206.6, 82.6)[l]{$e_1$\strut}}
  \put(515.89,407.67){\fontsize{11.76}{14.11}\selectfont \makebox(206.6, 82.6)[l]{$e_2$\strut}}
  \put(236.20,38.36){\fontsize{11.76}{14.11}\selectfont \makebox(124.0, 82.6)[l]{(a)\strut}}
  \put(251.36,158.25){\fontsize{11.76}{14.11}\selectfont \makebox(206.6, 82.6)[l]{$G_1$\strut}}
  \put(303.46,490.50){\fontsize{11.76}{14.11}\selectfont \makebox(206.6, 82.6)[l]{$x_3$\strut}}
  \put(111.29,435.88){\fontsize{11.76}{14.11}\selectfont \makebox(165.3, 82.6)[l]{$y_1$\strut}}
  \put(466.76,509.10){\fontsize{11.76}{14.11}\selectfont \makebox(206.6, 82.6)[l]{$e_3$\strut}}
  \put(1922.57,37.52){\fontsize{11.76}{14.11}\selectfont \makebox(124.0, 82.6)[l]{(c)\strut}}
  \put(117.36,324.61){\fontsize{11.76}{14.11}\selectfont \makebox(165.3, 82.6)[l]{$y_2$\strut}}
  \put(864.57,463.20){\fontsize{11.76}{14.11}\selectfont \makebox(206.6, 82.6)[l]{$e_5$\strut}}
  \put(1139.65,623.10){\fontsize{11.76}{14.11}\selectfont \makebox(206.6, 82.6)[l]{$e_1$\strut}}
  \put(1178.18,197.21){\fontsize{11.76}{14.11}\selectfont \makebox(206.6, 82.6)[l]{$e_3$\strut}}
  \put(813.41,323.77){\fontsize{11.76}{14.11}\selectfont \makebox(206.6, 82.6)[l]{$e_6$\strut}}
  \put(1184.49,427.86){\fontsize{11.76}{14.11}\selectfont \makebox(206.6, 82.6)[l]{$e_2$\strut}}
  \put(1113.53,51.30){\fontsize{11.76}{14.11}\selectfont \makebox(124.0, 82.6)[l]{(b)\strut}}
  \put(1031.43,553.00){\fontsize{11.76}{14.11}\selectfont \makebox(206.6, 82.6)[l]{$e_4$\strut}}
  \end{picture}%
  \caption{The graphs in the proof of Case 3 of Conjecture \ref{con1.5}.} \label{fig8}
\end{figure}

\noindent {\bf Subcase 3.1.} $c(G_2)=0$ and $c(G_1)=k-2$. Let
$|V(G_2)|=n_2$. Then we have $3n_2=\sum_{v\in
V(G_2)}d_{G}(v)=2(n_2-1)+3=2n_2+1$. Therefore, $n_2=1$, i.e.,
$G_2=S_1$. Let $V(G_2)=\{x\}$, $e_1=xx_1$, $e_2=xx_2$ and
$e_3=xx_3$. Let $N_{G_1}(x_2)=\{y_1, y_2\}$ (see Figure \ref{fig8}
(a)). Let $F'=\{e_1, e_3, x_2y_1, x_2y_2\}$. Then $G-F'=G'_1+ G'_2$,
where $G'_2\cong S_2$ and $G'_1$ is a graph obtained from $G_1$ by
deleting a vertex of degree $2$. Therefore, $c(G'_1)=k-3\geq 2$. If
$G'_1\not\cong K_{2,3}$, then we are done by Claim 1. If $G'_1\cong
K_{2,3}$, then $n=7$, which is a contradiction.

\noindent {\bf Subcase 3.2.} $c(G_2)=1$ and $c(G_1)=k-3$. Let
$|V(G_2)|=n_2$. Then we have $3n_2=\sum_{v\in
V(G_2)}d_{G}(v)=2n_2+3$. Therefore, $n_2=3$, i.e., $G_2$ is a
triangle. If $G_1\not\cong K_{2,3}$, then we are done by Claim 1. If
$G_1\cong K_{2,3}$, then $G$ must be the graph as given in Figure
\ref{fig8} (b). Let $F'=\{e_1, e_4, e_5, e_6\}$. Then $G-F'=G'_1+
G'_2$, where $G'_1\cong S_2$ and $G'_2$ is a bicyclic graph which is
not isomorphic to $K_{2,3}$. Then we are done by Claim 1.

\noindent {\bf Subcase 3.3.} $c(G_2)=2$ and $c(G_1)=k-4$. Let
$|V(G_2)|=n_2$. Then we have $3n_2=\sum_{v\in
V(G_2)}d_{G}(v)=2(n_2+1)+3=2n_2+5$. Therefore, $n_2=5$. If neither
$G_1$ nor $G_2$ is isomorphic to $K_{2,3}$, then we are done by
Claim 1. Otherwise, we assume that $G_2\cong K_{2,3}$ (similar for
$G_1\cong K_{2,3}$). Then $G$ must have the structure as given in
Figure \ref{fig8} (c). Let $F'=\{e_1, e_2, e_4, e_5\}$. Then
$G-F'=G'_1+ G'_2$, where $G'_2\cong K_{2,2}$ and $G'_1$ is a
$(k-4)$-cyclic graph which is not isomorphic to $K_{2,3}$. Then we
are done by Claim 1. The proof is thus complete. \qed


\begin{thebibliography}{s1}

\bibitem{BP}
R.B. Bapat, S. Pati, Energy of a graph is never an odd integer, {\it
Bull. Kerala Math. Assoc.\/} 1(2004), 129--132.

\bibitem{BM}
J.A. Bondy, U.S.R. Murty, {\it Graph Theory with Applications\/},
Macmillan London and Elsvier, New York (1976).

\bibitem{CDS}
D. Cvetkovi\'{c}, M. Doob, H. Sachs, {\it Spectra of Graphs --
Theory and Application\/}, Academic Press, New York, 1980.

\bibitem{CP}
D. Cvetkovi\'{c}, M. Petri\'{c}, A table of connected graphs on six
vertices, {\it Discrete Math.\/} 50 (1984), 37--49.

\bibitem{DS}
J. Day, W. So, Graph energy change due to edge deletion, {\it Lin.
Algebra Appl.\/} 428(2008), 2070--2078.

\bibitem{G1}
I. Gutman, On graphs whose energy exceeds the number of vertices,
{\it Lin. Algebra Appl.\/} 429(2008), 2670--2677.

\bibitem{GLSZ}
I. Gutman, X. Li, Y. Shi, J. Zhang, Hypoenergetic trees, {\it MATCH
Commun. Math. Comput. Chem.\/} 60(2009), 415--426.

\bibitem{GR}
I. Gutman, S. Radenkovi\'c, Hypoenergetic molecular graphs, {\it
Indian J. Chem.\/} 46A (2007), 1733--1736.

\bibitem{LM1}
X. Li, H. Ma, Hypoenergetic and strongly hypoenergetic $k$-cyclic
graphs, accepted for publication in {\it MATCH Commun. Math. Comput.
Chem.\/}

\bibitem{LM2}
X. Li, H. Ma, All hypoenergetic graphs with maximum degree at most
3, accepted for publication in {\it Lin. Algebra Appl.\/}

\bibitem{MKG}
S. Majstorovi\'{c}, A. Klobu\v{c}ar, I. Gutman, Selected topics from
the theory of graph energy: Hypoenergetic graphs, in: {\it
Applications of Graph Spectra\/}, Math. Inst., Belgrade, 2009,
65--105.

\bibitem{N}
V. Nikiforov, Graphs and matrices with maximal energy, {\it J. Math.
Anal. Appl.\/} 327(2007), 735--738.

\bibitem{N2}
V. Nikiforov, The energy of $C_4$-free graphs of bounded degree,
{\it Lin. Algebra Appl.\/} 428(2008), 2569--2573.

\end{thebibliography}
\end{document}